\documentclass[11pt,reqno]{article}
\usepackage{graphicx,amssymb}
\usepackage{amssymb,amsmath}
\usepackage{hyperref}
\usepackage{amsfonts}
\usepackage{amscd}
\usepackage{amsthm}
\usepackage{color}
\numberwithin{equation}{section}
\newtheorem{thm}{Theorem}[section]
\newtheorem{Lem}{Lemma}[section]

\newtheorem{Def}{Definition}[section]
\newtheorem{Rem}{Remark:}[section]
\usepackage{geometry}
\begin{document}
\title{Inverse coefficient problem for cascade system of fourth and second order partial differential equations}
\baselineskip 12.5pt
\author{\large Navaneetha Krishnan Murugesan\\ 
\small Department of Mathematics, Central University of Tamil Nadu, Thiruvarur 610005, INDIA\\ [2ex]
\large Sakthivel Kumarasamy \\ 
\small Department of Mathematics, Indian Institute of Space Science and Technology, Trivandrum 695547, INDIA\\[2ex]
\large Alemdar Hasanov \\ 
\small Department of Mathematics, Kocaeli University, Kocaeli 41001, TURKEY\\[2ex]
\large Barani Balan Natesan\footnote{Corresponding Author E-Mail: baranibalan@acad.cutn.ac.in}\\ 
\small Department of Mathematics, Central University of Tamil Nadu, Thiruvarur 610005, INDIA}
\date{\empty}
\maketitle
\begin{abstract}
The study of the paper mainly focusses on recovering the dissipative parameter in a cascade system coupling a bilaplacian operator to a heat equation from final time measured data via quasi-solution based optimization. The coefficient inverse problem is expressed as a minimization problem. We proved that minimizer exists and the necessary optimality condition which plays the crucial role to prove the required stability result for the corresponding coefficient is derived. Utilising the conjugate gradient approach, numerical results are examined to show the method's effectiveness.
\end{abstract}
\section{Introduction}
 Inverse and ill-posed problems have been studied in various fields of engineering and technology
during the last few decades (see \cite{arr}, \cite{hasan}, \cite{reddy} and references therein). Unlike others, inverse
coefficient problems have two distinguished features. First, these problems are severely ill-posed.
The second, these problems are nonlinear if even the differential operator in the direct problem
is linear. These properties make it difficult to study fundamental problems such as uniqueness
and stability. In this work, we are interested to study a coefficient inverse problem for coupled partial differential equations, one of second order and the other of fourth order, that form a parabolic system. This parabolic system combines dissipative and dispersive characteristics, represents front propagation in reaction-diffusion events, and simultaneously sustains a stable solitary pulse. The stabilised Kuramoto-Sivashinsky system was proposed in \cite{mal} and consists of a linearly coupled one-dimensional Kuramoto-Sivashinsky-Korteweg de Vries (KS-KdV) equation to an another dissipative equation, which has the form
\begin{eqnarray} \label{0.1}
 \hspace{-.2in}\left\{ \begin{array}{ll}
 u_{t}+ b u_{xxxx}+u_{xxx}+a u_{xx}+uu_x=v_x, \   & (x,t) \  \in \ \Omega_T,  \\[2mm]
 v_{t}- d v_{xx}+c v_x=u_x, \ & (x,t) \  \in \ \Omega_T.  
 \end{array}
 \right.\end{eqnarray}
where $b$ stands for long-wave instability, $a$ for the short-wave dissipation, and $d>0$ for the dissipative parameter (for instance, the loss constant) and $c$ denotes the group velocity mismatch between wave fields. This system serves as a one-dimensional model for the propagation of waves and turbulence.

 The current work studies the mathematical analysis of the inverse coefficient problem with final time measurement for the simplified linear version of \eqref{0.1} (see, \cite{cer})
  \begin{eqnarray} \label{dir}
 \hspace{-.2in}\left\{ \begin{array}{ll}
 u_{t}+ u_{xxxx}=v, \   & (x,t) \  \in \ \Omega_T,  \\[2mm]
 v_{t}- (d(x) v_{x})_{x}=0, \ & (x,t) \  \in \ \Omega_T,  \\[2mm]
 u(x,0)=u_0(x),\  v(x,0)=v_0(x), \ & x \in \Omega, \\[2mm]
 u(x,t)=  u_{xx}(x,t)=0,\ v(x,t)=0, \ & (x,t) \in \partial \Omega_T,
 \end{array}
 \right.\end{eqnarray}
where $u$ and $v$ represents two real wave fields and $d(x)>0$ denotes the dissipative parameter, $\Omega :=(0,1), \ \Omega_T:= \Omega\times(0,T), \ $ and $\partial \Omega_T := \partial\Omega \times (0,T)$. We need to reconstruct the dissipative parameter $d(x)$ from the following final time measured outputs
 \begin{eqnarray} \label{1.2}
 m_1(x):=u(x,T;d)\ \ \mbox{and}\ \ m_2(x):=v(x,T;d).
 \end{eqnarray}
 Here, we assume that $(m_1(x),m_2(x))$ is the measured output data that may contain noisy data and $u(x,T;d), \ v(x,T;d)$ are the output data corresponding to the given coefficient $d(x)$. The main objective of this study is to establish a stability estimate for the inverse coefficient problem governed by \eqref{dir} and \eqref{1.2}. We use the approach based on the weak solution theory combined with the adjoint method and the Tikhonov functional,proposed in \cite{ah} and develop this approach for the system \eqref{dir} governed by the cascade equations.
 
Let us have a quick review of literature related to this work. ICP/ISP related to parabolic equation have been widely investigated in the literature by various authors(see, \cite{Abdel}-\cite{les2},\cite{deng2}-\cite{erd},\cite{gnan}-\cite{hasa} and references therein). There are several works have been done in the case of inverse problem of system of PDEs and higher order PDEs. The authors in \cite{gnan} determined the source terms in Lotka-Volterra system with the knowledge of final time and Dirichlet measured data via quasi-solution approach. The retrieval of the dissipative parameter in a coupled parabolic-elliptic system from the final time measurement data was done in \cite{deng2}. With the help of aforementioned measurement data, the coefficients in reaction diffusion system have been identified in \cite{sak} and the coefficients along with initial data  have been identified in \cite{gnan}. The inverse coefficient problem for linear KdV equation from final time data was investigated in \cite{sak1} and the unknown coefficient of the same equation using Neumann boundary data also studied in \cite{sak2}. The authors \cite{hasa3},\cite{hasa2} determined the shear force acting on the boundary and time dependent source term from Euler-Bernoulli beam equation. Further, inverse problems for parabolic equation along with the numerical study can also be found in various articles \cite{Abdel}-\cite{Abdel1}, \cite{deng2}-\cite{erd} and references therein. In this paper, we apply the theory developed in \cite{ah}. According to this theory, a weak solution to the direct problem is first determined, since the input parameters in any inverse problem are usually non-smooth. Then the input-output operator is introduced and the properties of this operator are investigated using a priori estimates for the weak solution. Further, the Tikhonov functional is introduced and the inverse problem is reformulated as the problem of minimization of this functional. Here the key points are that, on the basis of this, the existence of a quasi-solution is proved, and formula for the gradient of the Tikhonov functional is introduced. The latter is an esstential tool for the numerical solution of the inverse problem. The idea behind the theory is to restrict the unknown parameter to some admissible set and transform the problem into a minimization problem by introducing a suitable auxiliary functional depending on an unknown function. Then we solve the minimization problem using the classical optimal control framework. This idea is widely used in the literature(see for instance \cite{les,les1,les2,chen,deng1,gnan1,hasa,hasa3,hasa1,
sak,sak2}).  
 
The boundary controllability of the system $\eqref{dir}$ is studied in \cite{cer} with the constant dissipative parameter when the control acts on the boundary of the heat equation. It was shown that the dissipative parameter plays a crucial role in determining the nature of the controllability results, namely, approximate controllability, null-controllability, failure of null-controllability, etc. By considering the importance of the dissipative parameter in the qualitative study of the system $\eqref{dir}$, in this article, we determine the dissipative parameter $d(x)$ from the final time measurements of $u \mbox{ and } v$. We also made a comparative study that shows how the stability constant varies with the change in the regularization parameter. Further, we have found an upper bound for the value of the final time under which Lipschitz type stability holds for the unknown coefficient $d(x)$. To the best of our knowledge, such a study has not been done in the literature for the cascade system $\eqref{dir}$. 

Let $H^m(\Omega),m\geq 0$ denotes the standard Sobolev space with norm denoted by $\|\cdot\|_{H^m(\Omega)}$. We employ the time-dependent function spaces $L^2(0,T;H^m(\Omega))$ which consists of all measurable functions from $(0,T)$ into $H^m(\Omega)$ such that their $H^m$-norm squared is integrable over $(0,T)$. Additionally, we use $L^{\infty}(0,T;H^m(\Omega))$ space as well. Moreover, unless specified $\|\cdot \|$ denotes $L^2(\Omega)$ norm.

We study the inverse problem with the following class of admissible dissipative parameters:
\begin{eqnarray}\label{1.3}
\lefteqn{\hspace{-2.5in}\mathcal{D}=\left\{ d\in H^1(\Omega): 0<\alpha_0\leq d(x)\leq \alpha_1, \ \forall x\in \Omega, \mbox{ and }\|d\|_{H^1(\Omega)}\leq \alpha_2\right\}.}
\end{eqnarray}

 Now, we formulate the inverse problem as a minimization problem. Consider an input-output operator, namely,  \\
\begin{eqnarray}\label{op}
\left\{ \begin{array}{ll} \Psi:\mathcal{D}\subset H^1(\Omega)\mapsto \mathcal{K} \\
 \Psi[d](x):=\left(u(x,t;d),v(x,t;d)\right)\left. \right|_{t=T},
\end{array} \right. 
\end{eqnarray}
where $\mathcal{K}=H^1_0(\Omega) \times L^{2}(\Omega)$. By using this operator, we reformulate the inverse coefficient problem(ICP) $\eqref{dir}-\eqref{1.2}$ into an operator equation  
  \begin{eqnarray}\label{op1}
 \Psi[d](x)=(m_1(x),m_2(x)),\  d\in \mathcal{D}, \ (m_1(x),m_2(x))\in L^2(\Omega) \times L^2(\Omega).
\end{eqnarray}

 Since the measured output always contains noise, an exact equality in $\eqref{op1}$ is impossible in practice. To this end we introduce the regularized Tikhonov functional
\begin{eqnarray}\label{3.1}
 {J_{\gamma}({d})}=\frac{1}{2}\int\limits_{\Omega}{(u(x,T;d)-m_1(x))^{2} dx}+\frac{1}{\textnormal{2}}\int\limits_{\Omega}{(v(x,T;d)-m_2(x))^{2} dx}+\frac{\gamma}{2} \|d\|_{H^1(\Omega)}^{2},
\end{eqnarray}
and study the inverse problem $\eqref{dir}-\eqref{1.2}$ as the following minimization problem
\begin{eqnarray}\label{3.1a}
J_{\gamma}({d^*})=\min\limits_{d\in \mathcal{D}} {J_{\gamma}({d})},
\end{eqnarray} 
where $(u(x,t;d),v(x,t;d)) \in L^2(0,T;V_0^2(\Omega))\times L^2(0,T;H_0^1(\Omega))$ is the weak solution of $\eqref{dir}$ which corresponds to the admissible class of dissipative parameters $\mathcal{D}$ defined in $\eqref{3.2}$. The term $\|d\|_{H^1(\Omega)}$ in $\eqref{3.1}$ represents the Tikhonov regularization term which stabilizes the resolution process of the output least square formulation for the corresponding parameter $\gamma>0$.

The following is how the paper is set up. Using the Faedo-Galerkin approach, we examine the well-posedness of the direct problem $\eqref{dir}$ and the associated adjoint problem  in section 2. In section 3, we prove the inverse problem $\eqref{dir}-\eqref{1.2}$ is illposed through the input-output operator and prove the existence of a minimizer. In section 4, we discuss the Fr\'echet derivative of the functional. In section 5, we derive the necessary optimality condition. The main result of the paper is established in the penultimate section in which we discuss the stability of the unknown coefficient  under suitable norms with the aid of necessary optimality condition, that is, we prove a small change/error in the given measurement does not cause much deviation in the unknown parameter $d(x).$ For this, we let $(u,v)$ and $(\widetilde{u},\widetilde{v})$ be the two solutions of the system \eqref{dir} for the corresponding coefficients $d(x) \mbox{ and } \widetilde{d}(x)$. Then for sufficiently small time $T,$ the following stability estimate holds 
$$\|d-\widetilde{d}\|_{H^1(\Omega)} \leq C\left(\|m_1-\widetilde{m}_1\|+\|m_2-\widetilde{m}_2\|\right),$$ 
where $C$ is a constant depending on $T$ and other parameters $m_1, m_2,\widetilde{m}_1,\widetilde{m}_2$. In the last section, we reconstruct the unknown coefficient of the inverse problem using Conjugate Gradient Method (CGM). We use Finite Difference scheme to solve the direct and adjoint systems and we compare the effects of Polak-Ribiere and Fletcher-Reeves conjugate coefficients in the numerical scheme. Algorithm is also provided and two numerical examples are presented in the corresponding section.
\section{Analysis of the Direct and Adjoint Problem}
   In this section, we investigate the existence and uniqueness result of the direct problem $\eqref{dir}$ and the corresponding adjoint problem.                                                                                                                                                                                                                                                 \begin{Def} [Weak solution]\label{D1.1}
Let $d \in H^1(\Omega),  u_0 \in L^2(\Omega), \mbox{ and } v_0 \in L^2(\Omega) $. A function $(u,v)\in L^2(0,T;V_0^2(\Omega))\times L^2(0,T;H_0^1(\Omega))$ with $(u_t,v_t)\in L^2(0,T;V_0^{2}(\Omega)^*)\times L^2(0,T;H^{-1}(\Omega))$ is referred to as a weak solution if $u(x,0)=u_0(x),v(x,0)=v_0(x)$ and the integral equality mentioned below holds 
\begin{eqnarray}\label{3.2}
\left\{\begin{array}{l}
\displaystyle{\langle u_{t},\phi \rangle + (u_{xx}, \phi_{xx})= (v,\phi),}\\ [4mm]
\displaystyle{\langle v_{t},\phi \rangle +(d(x) v_{x}, \phi_{x}) = 0, \  \forall \ \phi \in V_0^2(\Omega), \ \mbox{a.e.} \ t \in [0,T],}
\end{array}\right.
\end{eqnarray}
where $\phi \in V_0^2(\Omega):=\{z \in H^2(\Omega):z=0 \mbox{ on } \partial \Omega\}$ and $V_0^{2}(\Omega)^*$ is the dual space of $V_0^2(\Omega).$
\end{Def}
The verification of the initial data is justified below.
 \begin{thm}\label{thm1}
 Assume $u_0,v_0\in L^2(\Omega) \mbox{ and } d\in \mathcal{D}$. Then there exists a unique weak solution for $\eqref{dir}$ in the sense of Definition \ref{D1.1} satisfying the estimates :
 \begin{eqnarray}\label{2.1} 
 \|(u,v)\|^2_{X \times Y} &\leq& C_1(T,\alpha_0) (\|u_0\|^2+\|v_0\|^2),  \\
\label{2.1a} 
 \mbox{ and }\|(u_t,v_t)\|^2_{L^2(0,T;V_0^2(\Omega)^*) \times L^2(0,T;H^{-1}(\Omega))} &\leq&  C_2(T,\alpha_0,\alpha_1) \left(\|u_0\|^2 + \|v_0\|^2 \right),
 \end{eqnarray}
where $X=L^{\infty}(0,T;L^2(\Omega))\cap L^2(0,T;V_0^2(\Omega))$ and $Y=L^{\infty}(0,T;L^2(\Omega)) \cap L^2(0,T;H_0^1(\Omega)),$\\
$$C_1(T,\alpha_0)=\displaystyle \frac{(1+Te^T+4CTe^T)(1+\alpha_0)}{2\alpha_0}, \ C>0 \mbox{ is the Ehrling's Lemma constant and }$$ \\ $$C_2(T,\alpha_0,\alpha_1)=\displaystyle \left(1+3Te^T+\frac{\alpha_1^2(1+Te^T)}{2\alpha_0}\right) (\|u_0\|^2+\|v_0\|^2).$$
 \end{thm}
 \noindent{\bf Proof:} To prove the theorem, we apply the Faedo-Galerkin method to construct approximate solution in finite dimensional space and extend it to infinite dimensional space using standard limiting arguments.\\
 {\it Step $1$:} Let $(w_k)_{k\in \mathbb{N}}\ $ be an orthogonal basis in $V_0^2(\Omega)$ and an orthonormal basis of $L^2(\Omega)$. Now we construct the approximate solution spanned by the subspace $\{w_1,w_2, \ldots, w_m \}$ as follows
 \begin{eqnarray*}
 u_m(x,t)=\sum\limits_{k=1}^m{h_m^k(t)w_k(x)} \mbox{ and } 
 v_m(x,t)=\sum\limits_{k=1}^m{g_m^k(t)w_k(x)},
 \end{eqnarray*}
 where $h^k_m(t) \mbox{ and } g_m^k(t)$ are the solutions of the system of ordinary differential equations given by
 \begin{eqnarray} \label{2.2}
 \hspace{-.2in}\left\{
\begin{array}{ll}
 (u'_{m},w_k) +( u_{m,xx}, w_{k,xx})=(v_{m},w_k),  \\[2mm]
 (v'_{m},w_k) +(d(x) v_{m,x}, w_{k,x})=0,  \\[2mm]
 \ h_m^k(0)=(u_0,w_k),\ g_m^k(0)=(v_0,w_k),
 \end{array}
 \right.\end{eqnarray}
 for $k=1,2,\cdots,m,$ where $(\cdot,\cdot)$ denotes the inner product in $L^2(\Omega)$. By the standard theory of ODEs, we can prove that the ODEs associated with $\eqref{2.2}$ have a solution on $[0,T].$\\
 
\noindent  {\it Step $2$:} Multiply $\eqref{2.2}_a \mbox{ and } \eqref{2.2}_b$ (represent first and second equation of $\eqref{2.2}$) by $h_m^k(t) \mbox{ and } \\ g_m^k(t)$,$ \ k=1,2,\cdots m$ respectively, take sum from $1$ to $m$, we get
\begin{eqnarray*}
&&\frac{1}{2}\frac{d}{dt}\int\limits_{\Omega} |u_m|^2 dx+\int\limits_{\Omega} | u_{m,xx}|^2 dx = \int\limits_{\Omega} u_m v_m dx,\\
&&\frac{1}{2}\frac{d}{dt} \int\limits_{\Omega} | v_m|^2 dx+ \int\limits_{\Omega} d(x) | v_{m,x}|^2 dx =0.
\end{eqnarray*}
Applying Cauchy's inequality, assumption on the coefficient 'd' and adding them together, we get
\begin{eqnarray}\label{2.5}
 \frac{d}{dt} \left(\|u_m(t)\|^2 +\|v_m(t)\|^2 \right) + 2\| u_{m,xx}(t)\|^2 + 2\alpha_0 \| v_{m,x}(t)\|^2 \leq \|u_m(t)\|^2 +\|v_m(t)\|^2.
 \end{eqnarray}  
    Now applying Gronwall's inequality to \eqref{2.5}, we have
   \begin{eqnarray}
\label{2.6}  \|u_m(t)\|^2 +\|v_m(t)\|^2  \!\!&\leq &  \!\! e^{T} \left(\|u_0\|^2 +\|v_0\|^2 \right).
 \end{eqnarray}
   \noindent Integrating $\eqref{2.5}$ over $(0,t)$ and using $\eqref{2.6}$, to get
 \begin{eqnarray} 
\label{2.7}
2\int\limits_0^t\!\| u_{m,xx}(\tau)\|^2d\tau\!+ 2\alpha_0\int\limits_0^t\| v_{m,x}(\tau)\|^2d\tau \!&\!\!\leq \!\! & \left(1+Te^{T}\right)\left(\|u_0\|^2 +\|v_0\|^2 \right).
  \end{eqnarray}
By using Ehrling’s lemma (see Theorem 7.30 in \cite{ren}), for any $\epsilon \geq 0$, we obtain
\begin{eqnarray}\label{2.9 c}
\| u_m(t) \|^2_{H_0^1(\Omega)} \leq \epsilon (\| u_m(t) \|^2_{H_0^1(\Omega)}+\| u_{m,xx}(t) \|^2) + C(\epsilon) \| u_m(t) \|^2.
\end{eqnarray}
Integrating over $(0,t)$ and substituting \eqref{2.6} and \eqref{2.7} in \eqref{2.9 c} and choosing $\epsilon = \frac{1}{2}$, we get
\begin{eqnarray}\label{2.9 d}
\|u_m(t)\|^2_{L^2(0,T;H_0^1(\Omega))} \leq   \left( 2CT e^T+\frac{1+Te^T}{2} \right) (\|u_0\|^2+\|v_0\|^2).
\end{eqnarray}
Now fix any $z \in V_0^2(\Omega)$ with $\|z\|_{V_0^2(\Omega)} \leq 1$ and write $z=z^1 + z^2$ where $z^1 \in span  \{w_k\}_{k=1}^m$ and $(z^2,w_k)=0 \ (k=1,2,\cdots,m).$ Since the functions $\{w_k\}_{k=0}^{\infty}$ are orthogonal in $V_0^2(\Omega)$, $\|z^1\|_{V_0^2(\Omega)} \leq \|z\|_{V_0^2(\Omega)} \leq 1$, from \eqref{2.2}, we get
$$\langle u'_{m},z \rangle =(u'_{m},z^1)=(v_{m},  z^1)-(u_{m,xx},(z^1)_{xx})$$
Using H\"older's inequality, we get
$$  | \langle u'_{m},z \rangle | \leq \| v_{m}(t) \| \|z^1 \| + \|u_{m,xx}\| \|(z^1)_{xx}).$$ Squaring on both sides and then integrating over $(0,t)$ and using \eqref{2.6} and \eqref{2.7}, we get
\begin{eqnarray}\label{2.10a}
  \|u'_{m}(t)\|^2_{L^2(0,T;V_0^2(\Omega)^*)} \leq (1+3Te^{T})(\|u_0\|^2 + \|v_0\|^2),
  \end{eqnarray}
  since $\|z\|_{V_0^2(\Omega)} \leq 1$. Similarly, by taking functions $z\in H_0^1(\Omega)$ with $\|z\|_{H_0^1(\Omega)}\leq 1$, we can deduce the following estimate
 \begin{eqnarray}\label{2.10}
  \|v'_{m}(t)\|^2_{L^2(0,T;H^{-1}(\Omega))} \leq \frac{\alpha_1^2(1+Te^{T})}{2\alpha_0}(\|u_0\|^2 + \|v_0\|^2).
 \end{eqnarray}
Combining \eqref{2.6}-\eqref{2.7} and \eqref{2.9 d}-\eqref{2.10}, we get
\begin{eqnarray}\label{2.11}
&&  \|(u_m,v_m)\|^2_{X \times Y}\leq  C_1(T,\alpha_0) \left(\|u_0\|^2 + \|v_0\|^2\right)   \\\label{2.11a} && \|(u'_m,v'_m)\|^2_{L^2(0,T;V_0^2(\Omega)^*) \times L^2(0,T;H^{-1}(\Omega))} \leq  C_2(T,\alpha_0,\alpha_1) \left(\|u_0\|^2 + \|v_0\|^2\right).
\end{eqnarray}
{\it Step $3$:} From the above estimates, we conclude that
\begin{eqnarray*}
\left\{\begin{array}{l}
\{u_m\}_{m=1}^{\infty} \mbox{is bounded in } L^{\infty}(0,T;L^2(\Omega))\cap L^{2}(0,T;V_0^2(\Omega)),\\[2mm]
\{v_m\}_{m=1}^{\infty} \mbox{is bounded in } L^{\infty}(0,T; L^2(\Omega))\cap L^{2}(0,T;H^1_0(\Omega)),\\[2mm]
\{u'_{m}\}_{m=1}^{\infty} \mbox{is bounded in } L^{2}(0,T;V_0^2(\Omega)^*),\\[2mm]
\{v'_{m}\}_{m=1}^{\infty} \mbox{is bounded in } L^{2}(0,T;H^{-1}(\Omega)).
\end{array}
\right.
\end{eqnarray*}
As a consequence of the above, we can extract subsequences (which is denoted in the same way as the original sequences) such that 
\begin{eqnarray}\label{2.12a}
\left\{\begin{array}{l}
u_m \rightharpoonup u \mbox{ weakly* in } L^{\infty}(0,T;L^2(\Omega))\cap L^{2}(0,T;V_0^2(\Omega)),\\[2mm]
v_m\rightharpoonup v \mbox{ weakly* in } L^{\infty}(0,T; L^2(\Omega))\cap L^{2}(0,T;H^1_0(\Omega)),\\[2mm]
u'_{m} \rightharpoonup u_t \mbox{ weakly in } L^{2}(0,T;V_0^2(\Omega)^*),\\[2mm]
v'_{m}\rightharpoonup v_t \mbox{ weakly in } L^{2}(0,T;H^{-1}(\Omega)), \mbox{ as } m\rightarrow \infty. 
\end{array}
\right.
\end{eqnarray}
 Now we fix an integer $N$ and choose a function ${\bf y}\in C^1([0,T];V^2_0(\Omega))$ defined by
\begin{eqnarray}\label{2.12}
{\bf y}(t)=\sum\limits_{k=1}^N l^k(t)w_k,
\end{eqnarray}
where $\{l^{k}(t)\}_{k=1}^N$ are given smooth functions. For $m\geq N$, multiply $\eqref{2.2}$ by $l^k(t)$, take summation over $k$, integration with respect to $t$ over $(0,T)$ and applying \eqref{2.12a}
 we obtain that $(u,v)$ solves \eqref{3.2} for all $y\in C^1([0,T];V_0^2(\Omega))$.\\
Since, the functions in $C^1([0,T];V_0^2(\Omega))$ are dense in $L^2(0,T;V^2_0(\Omega))$ the identities \eqref{3.2} holds for every ${\bf y }\in L^2(0,T;V^2_0(\Omega))$. 
 Since, $(u,v)\in L^{2}(0,T;V_0^2(\Omega))\times L^{2}(0,T;H_0^1(\Omega))$ and $(u_t,v_t)\in L^{2}(0,T;V_0^2(\Omega)^*)\times L^{2}(0,T;H^{-1}(\Omega)) $, we have $(u,v)\in C([0,T];H_0^{1}(\Omega))\times C([0,T];L^2(\Omega))$. We can verify the initial conditions $u(x,0)=u_0, \ v(x,0)=v_0$ in a standard way.
\\To prove uniqueness, we take two weak solutions, say, $(u_1,v_1)$ and $(u_2,v_2)$ of \eqref{3.2}. By taking $(u,v)=(u_1-u_2,v_1-v_2) $ which satisfies $\eqref{dir}$ with homogeneous initial and boundary conditions . With the help of energy estimates \eqref{2.6}, the uniqueness result can be proved. Hence the proof.\hfill{$\Box$}\\
Next, we prove the regularity of weak solutions of \eqref{dir}.
\begin{Lem}\label{lem.a}
Assume $(u_0,v_0)\in H^2(\Omega) \times H_0^1(\Omega) \mbox{ and } d \in \mathcal{D}$. Then $u \in C([0,T];H_0^1(\Omega))$ and $v \in C([0,T];L^2(\Omega))$.
\end{Lem}
\noindent{\bf Proof:}
Multiply $\eqref{dir}_a$ by $u_t$, integrate over $\Omega$, apply integration by parts several times and using Cauchy's inequality, we get 
 \begin{eqnarray}\label{2.8}
\|u_{t}(t)\|^2 + \frac{d}{dt} \| u_{xx}(t) \|^2 \leq \|v(t)\|^2
\end{eqnarray}
Integrating \eqref{2.8} over $(0,t)$ and using \eqref{2.6}, we get 
\begin{eqnarray}\label{2.9}
 \int\limits_0^t \|u_{\tau}(\tau)\|^2 d \tau + \| u_{xx} (t)\|^2 \leq T e^{T} (\|u_0\|^2 + \|v_0\|^2) + \| (u_0)_{xx} \|^2.
 \end{eqnarray}
  Multiply $\eqref{dir}_b$ by $v_t$, integrate over $\Omega$, and apply integration by parts, we obtain 
 \begin{eqnarray}\label{2.9 a}
\|v_{t}(t)\|^2 +\frac{\alpha_0}{2} \frac{d}{dt} \|v_{x}(t) \|^2 \leq \| v_{x}(t) \|^2
\end{eqnarray}
 Now applying Gronwall's inequality to \eqref{2.9 a}, we have
\begin{eqnarray} \label{2.9 b} 
\| v_{x}(t) \|^2 \leq e^{\frac{2T}{\alpha_0}}\| (v_0)_{x} \|^2.
\end{eqnarray}
Using Theorem \ref{thm1} and above estimates, we can arrive at the required result of Lemma \ref{lem.a}.       \hfill{$\Box$} \\
Next we study the adjoint system of \eqref{dir}. Let $(p,q)$ be the weak solution of the following system
\begin{eqnarray} \label{adj}
 \hspace{-.2in}\left\{ \begin{array}{ll}
 -p_{t}+ p_{xxxx}=0, \ (x,t) \  \in \ \Omega_T,  \\[2mm]
 -q_{t}- (d(x) q_{x})_{x}=p, \ (x,t) \  \in \ \Omega_T,  \\[2mm]
 p(x,T)=p_T(x),\  q(x,T)=q_T(x), \ x \in \Omega, \\[2mm]
 p(x,t)= p_{xx}(x,t)=0,\ q(x,t)=0, \ (x,t) \in \partial\Omega_T,
 \end{array}
 \right.\end{eqnarray}
where $p_T \in L^2(\Omega), \ q_T\in L^2(\Omega)$ are arbitrary functions.
 We can straight away state the well posedness of the adjoint problem (final value problem) by converting it to an appropriate initial value problem by setting $p(x,t;d)=\varphi(x,T-t;d), \ q(x,t;d)=\zeta(x,T-t;d)$ which satisfies the following initial boundary value problem
 \begin{eqnarray} \label{adj1}
 \hspace{-.2in}\left\{ \begin{array}{ll}
 \varphi_{t}+ \varphi_{xxxx} =0, \ (x,t) \  \in \ \Omega_T,  \\[2mm]
 \zeta_{t}- (d(x) \zeta_{x})_{x}=\varphi, \ (x,t) \  \in \ \Omega_T,  \\[2mm]
 \varphi(x,0)=p_T(x),\  \zeta(x,0)=q_T(x), \ x \in \Omega, \\[2mm]
 \varphi(x,t)= \varphi_{xx}(x,t)=0,\ \zeta(x,t)=0, \ (x,t) \in \partial \Omega_T.
 \end{array}
 \right.\end{eqnarray}
 The well-posedness of \eqref{adj} (or) \eqref{adj1} can be completed by the similar arguments of \\Theorem \ref{thm1}.

Now, we are going to derive a priori estimate for \eqref{adj} which will be used later.
\begin{Lem}\label{lem.b}
Let $(m_1,m_2) \in \mathcal{K}$ and $(p(x,t),q(x,t))$ be the solution of $(\ref{adj})$. Then the inequalities given below holds:
\begin{eqnarray}\label{4.8}
&& \max\limits_{t\in [0,t)}\left(\|p(t)\|^2+\|q(t)\|^2 \right) \leq e^{T}\left( \|p_T\|^2+\|q_T\|^2\right),\\ \label{4.8a}
&& \|p_{xx}\|^2_{L^2(0,T;L^2(\Omega))}+\alpha_0\|q_x\|^2_{L^2(0,T;L^2(\Omega))} \leq \left(\frac{1+Te^{T}}{2}\right)\left(\|p_T\|^2+\|q_T\|^2 \right).
\end{eqnarray}
\end{Lem}
\noindent{\bf Proof:}
By following the same procedure as in Theorem \ref{thm1}, we can obtain \eqref{4.8} and \eqref{4.8a}.
\hfill{$\Box$} 
 \section{Existence of Solutions for the Inverse Problem}
In this section, we prove that the operator $\Psi$ is a compact operator and Lipschitz continuous. Then we establish the existence of minimizer for the functional $J_\gamma$.

 In the case of noisy free data, ICP $\eqref{dir}-\eqref{1.2}$ can be reduced to $\eqref{op1}$. The inverse problem  is ill-posed, if the operator $\Psi$ is compact which is proved in the following lemma.
\begin{Lem}\label{lem1}
Let the conditions in Lemma $\eqref{lem.a}$ hold. Then the input-output operator $\Psi$ defined in $\eqref{op}$ is compact from $\mathcal{K}$ to $L^2(\Omega) \times L^2(\Omega)$.
\end{Lem}
\noindent{\bf Proof:} To prove $\Psi$ is compact, it is enough to show that the sequence of output data $\{u(x,T;d_m),\\v(x,T;d_m)\}$ is bounded in $\mathcal{K}$ for every bounded sequence $\{d_m\}\subset H^1(\Omega)$ since $\mathcal{K}\subset L^2(\Omega)\times L^2(\Omega)$ is compact. From $\eqref{1.3}$, the sequence $\{d_m\}$ is bounded in $\mathcal{D}$. Let  $\{u(x,T;d_m),v(x,T;d_m)\}$ denotes the corresponding output solution sequence of the direct problem $\eqref{dir}$. The estimates \eqref{2.1} and \eqref{2.9} show that
$$\max\limits_{t \in [0,T]} \|u(t)\|_{H_0^1(\Omega)} \leq C(\|u\|_{L^2(0,T;H^2(\Omega))}+\|u_t\|_{L^2(0,T;L^2(\Omega))}) \leq C, \ \max\limits_{t \in [0,T]} \|v(t)\| \leq C. $$
The above estimates imply that the sequence $\{u(x,T;d_m),v(x,T;d_m)\}$ is uniformly bounded in $\mathcal{K}$. This proves that $\Psi$ is a compact operator.\hfill{$\Box$}
\begin{Lem}\label{lem2}
Assume $u_0\in H^2(\Omega), \  v_0\in H^1_0(\Omega)$ and $d\in \mathcal{D}$. Then the input-output operator defined in $\eqref{op}$ is Lipschitz continuous, that is,
$$\|\Psi[d_1]-\Psi[d_2]\|\leq L_0\|\delta d\|_{H^1(\Omega)},$$
where $L_0=\displaystyle \left(\frac{e^T(1+Te^T)}{\alpha_0^2}(\|u_0\|^2+\|v_0\|^2)\right)^{\frac{1}{2}},$ \\[2mm]
 $\|\Psi[d_1]-\Psi[d_2]\|=\|(\delta u(\cdot,T),\delta v(\cdot,T))\|=\|\delta u(\cdot,T)\|+\|\delta v(\cdot,T)\|.$
\end{Lem}
{\bf Proof:} Let $(u_1(x,t;d_1),v_1(x,t;d_1)), (u_2(x,t;d_2),v_2(x,t;d_2))$ be the solutions of $\eqref{dir}$. Denote $\delta u:=u_1(x,t;d_1)-u_2(x,t;d_2), \delta v:=v_1(x,t;d_1)-v_2(x,t;d_2),$ and $\delta d=d_1-d_2$ which satisfies the following initial boundary value problem
 \begin{eqnarray} \label{4.1}
 \left\{ \begin{array}{ll}
\delta u_{t}+ \delta u_{xxxx}=\delta v, \  (x,t) \in \Omega_T,  \\[2mm]
\delta v_{t} - (d_2 \delta v_{x})_{x} = (\delta d (v_1)_x)_{x} , \  (x,t) \in \Omega_T,  \\[2mm]
 \delta u(x,0)=0,\ \delta v(x,0)=0, \ x\in \Omega, \\[2mm]
 \delta u(x,t)=0, \delta u_{xx}(x,t)=0, \ \delta v(x,t)=0,  \ (x, t) \in \partial \Omega_T.
 \end{array}
 \right.\end{eqnarray}
  Multiply $(\ref{4.1})_a \mbox{ and } (\ref{4.1})_b$ by $\delta u \mbox{ and } \delta v$ respectively, integrate over $\Omega$ and use integration by parts several times to arrive at
 \begin{eqnarray}
 \frac{1}{2}\frac{d}{dt}\!\int\limits_{\Omega}\!\!\left(|\delta u|^2\! +|\delta v|^2 \right)dx+ \!\int\limits_{\Omega}\! | \delta u_{xx}|^2 dx+\!\int\limits_{\Omega}\!\! d_2(x) | \delta v_{x}|^2 dx\!=\!\!\int_{\Omega} \delta u \delta v dx -\! \int_{\Omega} \delta d (v_1)_{x} \delta v_{x} dx.
 \end{eqnarray}
Now applying Cauchy's inequality and H\"older's inequality, we obtain
\begin{eqnarray}\label{4.2}
\frac{d}{dt}\left(\|\delta u(t)\|^2 +\|\delta v(t)\|^2 \right)+ 2\| \delta u_{xx}(t)\|^2 +\alpha_0\| \delta v_{x}(t)\|^2 \! \nonumber\\
  &&\hspace{-2.3in}\leq \left(\|\delta u(t)\|^2 +\|\delta v(t)\|^2  \right)+\frac{1}{\alpha_0} \|\delta d\|^2_{L^\infty(\Omega)} \| (v_1)_x(t) \|^2. 
  \end{eqnarray}
Using Gronwall's inequality to \eqref{4.2} and substituting \eqref{2.7} yield
\begin{eqnarray}\label{4.3}
 \max\limits_{t \in [0,T]} \left(\|\delta u(t)\|^2 +\|\delta v(t)\|^2 \right) \leq \frac{e^T(1+Te^T)}{2\alpha_0^2}\left(\| u_0\|^2 +\|v_0\|^2 \right) \|\delta d\|^2_{L^\infty(\Omega)}.
 \end{eqnarray}
 Integrating \eqref{4.2} over $(0,T)$ and using \eqref{4.3}, we get
 \begin{eqnarray}\label{3.9}
 2\int\limits_0^T \| \delta u_{xx} (\tau)\|^2 d\tau + \alpha_0 \int\limits_0^T \| \delta v_{x}(\tau)\|^2 d\tau \leq \frac{(1+Te^T)^2}{2\alpha_0^2} (\| u_0 \|^2 + \| v_0 \|^2)\| \delta d \|^2_{L^{\infty}(\Omega)}.
 \end{eqnarray}
From estimate $\eqref{4.3}$, we have
\begin{eqnarray}\label{3.8}
\|\Psi[d_1]-\Psi[d_2]\| \leq  \left(\frac{e^T(1+Te^T)}{2\alpha_0^2}(\|u_0\|^2+\|v_0\|^2)\right)^{\frac{1}{2}} \|\delta d\|_{L^{\infty}(\Omega)}.
 \end{eqnarray}
 By employing the embedding $H^1(\Omega)\hookrightarrow L^{\infty}(\Omega)$  to RHS of $\eqref{3.8}$, we attain 
$$ \|\Psi[d_1]-\Psi[d_2]\| \leq  L_0 \|\delta d\|_{H^1(\Omega)},$$
 which is the reqiured result. \hfill{$\Box$}\\
The following theorem ensures the existence of minimizer for the problem over the admissible set $\mathcal{D}$.
\begin{thm}
 Suppose the assumptions of Theorem \ref{thm1} hold. Then there exists a $d^* \in \mathcal{D}$ solving the minimization problem $\eqref{3.1a}$.
\end{thm}
\noindent{\bf Proof:} From $\eqref{3.1}$, we can see that $J_{\gamma}(d)$ is bounded below and $\displaystyle \frac{\gamma}{2} \|d\|_{H^1(\Omega)}^2$ for any $d \in  \mathcal{D}$\ which acts as a greatest lower bound for the functional $\eqref{3.1}$ since, the first term in RHS of $\eqref{3.1}$ is non-negative. Then from the definition of infimum, we can say that there exists a minimizing sequence $\{d_k\}\in \mathcal{D} \mbox{ of } J_{\gamma}$ such that 
\begin{eqnarray} \label{22}
\lim\limits_{k\rightarrow\infty}J_{\gamma}(d_k) = \inf\limits_{d \in \mathcal{D}} J_{\gamma}(d).
\end{eqnarray}
From the definition of  $\mathcal{D}$, we have $\|d\|_{H^1(\Omega)}\leq \alpha_2 $ from which we can extract subsequence (which is denoted as $\{d_k\}$) such that $d_k \rightharpoonup  d^*$ in $H^1(\Omega)$. The limit $d^*\in \mathcal{D}$, since the admissible set of dissipative parameters $\mathcal{D}$ is a closed convex subset of Hilbert space, so it is weakly closed. Further, from Theorem \ref{thm1} we can conclude that there exists a subsequence of $\{u(x,t;d_k)\}, \ \{v(x,t;d_k)\}$ such that $(u(x,t;d_k),v(x,t;d_k))\rightharpoonup (u^*(x,t),v^*(x,t))$ weakly in $L^2\left(0,T;V_0^2(\Omega) \right)\times L^2\left(0,T;H^1_0(\Omega)\right)$.
\par We next prove that $(u^*(x,t),v^*(x,t))= (u(x,t;d^*), v(x,t;d^*))$. For this, consider \\ $z \in C^1[0,T]$ such that $z(T)=0.$ Multiply $\eqref{3.2}$ by $z$ and integrate over $(0,T)$ to get
\begin{eqnarray}
\label{3.3}\lefteqn{\hspace{-3in}-\int\limits_{0}^T\!\!\int\limits_{\Omega}u(d_k)z_t\phi dx dt+\int\limits_{0}^T\int\limits_{\Omega} z u_{xx}(d_k) \phi_{xx} dx dt - \int\limits_{0}^T\!\!\int\limits_{\Omega}v(d_k)z \phi dx dt=\int\limits_{\Omega}u_0z(0)\phi dx,}\\
\label{3.4}\lefteqn{\hspace{-3in}-\int\limits_{0}^T\!\!\int\limits_{\Omega}\!v(d_k)z_t\phi dx dt\!+\!\int\limits_{0}^T\!\!\int\limits_{\Omega}\!(d_k\!-\!d^*) z \phi_{x} v_{x}(d_k) dx dt +\!\!\! \int\limits_{0}^T\!\!\int\limits_{\Omega}\! d^* z \phi_{x} v_{x}(d_k) dx dt \! =\!\! \int\limits_{\Omega} v_0 z(0)\phi dx,}\\
\lefteqn{\hspace{-2.9in} \forall \phi \in V_0^2(\Omega).} \nonumber 
\end{eqnarray}
The convergence of the second term on the L.H.S of $\eqref{3.4}$ is obtained as follows :
\begin{eqnarray*}
\int\limits_{0}^T\!\!\int\limits_{\Omega}(d_k-d^*)z \phi_{x} v_{x}(d_k) dx dt\!\!\!\!&\leq &\!\!\! \!\ T\| d_k-d^*\|_{L^2(\Omega)}\|z\|_{L^{\infty}(0,T)}\| \phi_{x}\|_{L^{\infty}(\Omega)}\|v(d_k)\|_{L^{\infty}(0,T;H_0^1(\Omega))}\\
&\rightarrow & 0 \mbox{ as } k\rightarrow \infty,
\end{eqnarray*}
since $d_k \rightarrow d$ in $L^2(\Omega)$ and $v(d_k)\in L^{\infty}(0,T;H_0^1(\Omega))$. If we allow $k\rightarrow \infty$ in $\eqref{3.3},\eqref{3.4}$, we get
\begin{eqnarray}\label{3.4a}
-\int\limits_{0}^T\int\limits_{\Omega}u^* z_t \phi dx dt+\int\limits_{0}^T\int\limits_{\Omega} u_{xx}^* z \phi_{xx} dx dt - \int\limits_{0}^T\int\limits_{\Omega}v^* z \phi dx dt&=&\int\limits_{\Omega}u_0z(0)\phi dx,\\\label{3.4b}
-\int\limits_{0}^T\int\limits_{\Omega}v^* z_t \phi dx dt+ \int\limits_{0}^T\int\limits_{\Omega} d^{*} z \phi_{x} v_{x}^* dxdt &=& \int\limits_{\Omega} v_0 z(0)\phi dx.  
\end{eqnarray}
On the other hand, consider \eqref{3.2} satisfied by $(u^*,v^*)$, time integration by parts and comparing with \eqref{3.4a} and \eqref{3.4b}, we get $u^*(x,0)=u_0(x),v^*(x,0)=v_0(x)$. From these we can conclude that $u^*(x,t)=u(x,t;d^*), v^*(x,t)=v(x,t;d^*).$ Finally, we show that the optimal for the functional $\eqref{3.1}$ is in fact $d^*(x)$.\\ 
It is clear that
\begin{eqnarray*}
\lefteqn{\hspace{-2.5in} \int\limits_{\Omega}\{(u(x, T;d_k) -m_1(x))-(u(x, T;d^*) -m_1(x))\}^2dx}\\
\lefteqn{\hspace{-1.5in}+\int\limits_{\Omega}\{(v(x, T;d_k) -m_2(x))-(v(x, T;d^*) -m_2(x))\}^2dx \geq 0.}
\end{eqnarray*}
As we know that $u(x,T;d_k)\rightharpoonup u(x,T;d^*),v(x,T;d_k)\rightharpoonup v(x,T;d^*)$ in $L^2(\Omega)$, we obtain

\begin{eqnarray}
\lefteqn{ \hspace{-3.09in}\liminf\limits_{k\rightarrow\infty}\int\limits_{\Omega}\left\{(u(x, T;d_k) -m_1(x))^2 + (v(x, T;d_k) -m_2(x))^2\right\}dx} \nonumber\\
\lefteqn{ \hspace{-2.8in}\geq 2 \lim\limits_{k\rightarrow\infty}\int\limits_{\Omega}(u(x, T;d_k)\! -\!m_1(x))(u(x, T;d^*)\! -\!m_1(x))dx -\int\limits_{\Omega}(u(x, T;d^*) \!-\!m_1(x))^2 dx}\nonumber\\
\lefteqn{\hspace{-2.53in}+2\lim\limits_{k\rightarrow\infty}\!\int\limits_{\Omega}\!(v(x, T;d_k) \!-\!m_2(x))(v(x, T;d^*)\! -\!m_2(x)) dx-\!\int\limits_{\Omega}\!(v(x, T;d^*) \!-\!m_2(x))^2 dx} \nonumber\\
\lefteqn{\hspace{-2.8in} = \int\limits_{\Omega}(u(x, T;d^*) -m_1(x))^2 dx+\int\limits_{\Omega}(v(x, T;d^*)-m_2(x))^2 dx.}\nonumber
\end{eqnarray}
The weak convergence of $d_k$ in $H^1(\Omega)$ and the lower semi-continuity in $L^2$ norm provides
\begin{eqnarray}
\liminf\limits_{k\rightarrow\infty} J_{\gamma}(d_k)\!\!\!\!&=&\!\!\!\!\liminf\limits_{k\rightarrow\infty}\!\!\left[\frac{1}{2}\!\int\limits_{\Omega}\!\!\{(u(x,\!T\!;d_k)\!-m_1(x))^2 \!+\!\! (v(x,\!T\!;d_k)\!-m_2(x))^2\}dx\!+ \frac{\gamma}{2} \|d_k\|^2_{H^1(\Omega)}\!\right]\nonumber\\
\!\!\!\!&\geq & \!\!\!\!\frac{1}{2}\int\limits_{\Omega}(u(x,T;d^*)-m_1(x))^2 dx+\frac{1}{2}\int\limits_{\Omega}(v(x,T;d^*)-m_2(x))^2 dx + \frac{\gamma}{2}\|{d^*}\|^2_{H^1(\Omega)}\nonumber\\
\!\!\!\!&=&\!\!\!\!J_{\gamma}(d^*).\nonumber
\end{eqnarray}
Using the last inequality and (\ref{22}), we get
$$\inf\limits_{d\in \mathcal{D}}J_{\gamma}(d)\leq J_{\gamma}(d^*)\leq \liminf\limits_{k\rightarrow \infty}J_{\gamma}(d_k)=\lim\limits_{k\rightarrow\infty}J_{\gamma}(d_k)=\inf\limits_{d \in \mathcal{D}}J_{\gamma}(d).$$
Hence $d^*$ is the minimizer of the objective functional $J_\gamma (d)$ in the admissible dissipative parameters $\mathcal{D}$ defined above in $\eqref{1.3}$. This concludes the proof. \hfill{$\Box$}
\section{Gradient of the Functional}
  \quad This section deals with the derivation of the Fr\'echet gradient of the objective functional $J_{\gamma}(d)$. 

 Assume the coefficients $d,d+\delta d \in \mathcal{D}$ and consider the increment in the functional \eqref{3.1} as follows:
 \begin{eqnarray}\label{4.9}
\lefteqn{\hspace{-3.0in} \delta {J}_\gamma(d):={J_\gamma}(d+\delta d)-{J_\gamma}(d) } \nonumber \\
 \lefteqn{\hspace{-2.446in}=\int\limits_{\Omega}\!\left(u(x,T;d)-m_1(x)\right)\delta  u(x,T) dx+\int\limits_{\Omega}\!\left(v(x,T;d)-m_2(x)\right)\delta  v(x,T) dx}  \nonumber \\ 
\lefteqn{\hspace{-2.0in} \!\!+ \gamma \Big(d,\delta d\Big)_{H^1(\Omega)} \!+\!\frac{1}{2}\int\limits_{\Omega}|\delta u(x,T)|^2 dx +\!\frac{1}{2}\int\limits_{\Omega}|\delta v(x,T)|^2 dx + \frac{\gamma}{2}\| \delta d\|^2_{H^1(\Omega)},}
 \end{eqnarray} 
where $(\delta u(x,T), \delta v(x,T))$ are the solutions of $\eqref{4.1}$ and 
$$\displaystyle  \Big(d,\delta d\Big)_{H^1(\Omega)}=\int\limits_{\Omega} d(x) \delta d(x) dx+\int\limits_{\Omega} d^{\prime}(x) \delta d^{\prime}(x) dx $$
 is the inner product defined over $H^1(\Omega).$ \\
By the formal Lagrangian method, the final time data of \eqref{adj} can be obtained as \\$p(x,T)=u(x,T;d)-m_1(x)$ and $q(x,T)=v(x,T;d)-m_2(x)$. Therefore, we have
\begin{eqnarray}
&&\hspace{-0.6in}\int\limits_{\Omega}\!\left(u(x,T;d)-m_1(x)\right)\delta  u(x,T) dx+\int\limits_{\Omega}\!\left(v(x,T;d)-m_2(x)\right)\delta  v(x,T) dx \nonumber\\
&=&\!\!\int\limits_{\Omega}\!p(x,T) \delta  u(x,T) dx+\int\limits_{\Omega}\!q(x,T)\delta  v(x,T) dx := I.\nonumber 
\end{eqnarray}
Multiplying the direct problem $\eqref{4.1}$ by $(p,q)$ and adjoint problem $\eqref{adj}$ by $(\delta u, \delta v)$ and then applying integration by parts several times, by using the initial and boundary conditions, we get
\begin{eqnarray}
I&=&\!\!\int\limits_{0}^T\!\!\int\limits_{\Omega}\Big( p_t(x,t) \delta  u(x,t) +p(x,t) \delta  u_t(x,t)+q_t(x,t)\delta  v(x,t)+q(x,t)\delta  v_t(x,t)\Big) dx dt\nonumber \\
\label{4.10}&=&-\int\limits_{0}^T\!\!\int\limits_{\Omega}\Big( \delta v_{x} q_{x} \delta d(x)+ q_{x} v_{x} \delta d(x) \Big)dx dt.
\end{eqnarray}
\begin{thm}\label{T4.1}
Assume that Dirichlet measured output $(m_1,m_2) \in L^2(\Omega) \times L^2(\Omega)$. Then the regularized Tikhonov functional $J_{\gamma}(d)$ is Fr\'echet differentiable on the set of admissible coefficients $\mathcal{D}$. The Frechet derivative at $d \in \mathcal{D}$ is given by the equation 
\begin{eqnarray}\label{4.15}
J^{\prime}_{\gamma}(d)=\Lambda + \gamma d,
\end{eqnarray}
 where $\Lambda$ is the solution of the following ordinary differential equation 
\begin{eqnarray}\label{4.14}
\left\{\begin{array}{l}
\displaystyle{\Lambda''-\Lambda=\int_0^T \!\! q_{x} v_{x}dt,} \\[3mm]
\displaystyle{\Lambda'(0)=\Lambda'(1)=0.}
\end{array}\right.
\end{eqnarray}
\end{thm}
\noindent{\bf Proof:} To prove the theorem, substitute $\eqref{4.10}$ in $\eqref{4.9}$ to get the following equality
\begin{eqnarray}\label{4.11}
\delta{J}_\gamma(d)\!\!&=&-\int\limits_{0}^T\!\!\int\limits_{\Omega} v_{x} q_{x}\delta d(x) dx dt -\int\limits_{0}^T\!\!\int\limits_{\Omega} \delta v_{x} q_{x}\delta d(x) dx dt+ \gamma \Big(d,\delta d\Big)_{H^1(\Omega)} \nonumber  \\  
\lefteqn{\hspace{0.4in}+\frac{1}{2}\int\limits_{\Omega}|\delta u(x,T)|^2 dx +\frac{1}{2}\int\limits_{\Omega} |\delta v(x,T)|^2 dx + \frac{\gamma}{2}\| \delta d\|_{H^1(\Omega)}^2.}
\end{eqnarray}
 By using $\eqref{4.3},\eqref{3.9},\eqref{4.8a}$ and applying Cauchy's inequality and H\"older's inequality, second, fourth and fifth terms on the RHS of \eqref{4.11} can be written as follows:
\begin{eqnarray}\label{4.15a}
\lefteqn{\hspace{-2.8in}-\int\limits_{0}^T\!\!\int\limits_{\Omega} \delta v_{x}  q_{x}\delta d(x) dx dt+\frac{1}{2}\int\limits_{\Omega}|\delta u(x,T)|^2 dx+ \frac{1}{2}\int\limits_{\Omega}\!|\delta v(x,T)|^2 dx}\nonumber\\ 
\lefteqn{\hspace{-2.9in}\leq \frac{1}{2}\int\limits_0^T\|\delta d\|^2_{L^{\infty}(\Omega)}\| q_{x}(\tau)\|^2 d \tau+\frac{1}{2}\int\limits_0^T\| \delta v_{x}(\tau) \|^2 d \tau+\frac{1}{2}\|\delta u(x,T)\|^2 + \frac{1}{2}\!\|\delta v(x,T)\|^2 }\nonumber \\
\lefteqn{\hspace{-2.9in}\leq \displaystyle \left(\frac{1+e^T(2\alpha_0 +T)+\alpha_0^2}{4\alpha_0^3}\right)(1+Te^T)\left(\|p_T\|^2+\|q_T\|^2+\|u_0\|^2+\|v_0\|^2 \right)\|\delta d\|^2_{L^{\infty}(\Omega)}. }
\end{eqnarray}
From $\eqref{2.6}$, we have
$\|p_T\|^2+\|q_T\|^2 \leq 2\left(2e^T\|u_{0}\|^2+2e^T\|v_{0}\|^2+\|m_1\|^2+\|m_2\|^2\right)$. Substitute it in $\eqref{4.15a}$ and apply the embedding result $H^1(\Omega)\hookrightarrow L^{\infty}(\Omega)$, we get
\begin{eqnarray} \label{4.13a}
\lefteqn{\hspace{-3.0in}-\int\limits_0^T \int\limits_{\Omega}\! \delta v_{x} q_{x} \delta d(x) dx dt\!+\!\frac{1}{2}\int\limits_{\Omega}|\delta u(x,T)|^2 dx \!+\! \frac{1}{2}\int\limits_{\Omega}|\delta v(x,T)|^2 dx \leq C_3\|\delta d\|^2_{H^1(\Omega)},}
\end{eqnarray}
where $C_3=\displaystyle \left(\frac{1+e^T(2\alpha_0 +T)+\alpha_0^2}{\alpha_0^3}\right)(1+Te^T)\left((2e^T+1)(\|u_0\|^2+\|v_0\|^2)+\|m_1\|^2+\|m_2\|^2 \right)$. Substitute \eqref{4.13a} in \eqref{4.11} to arrive at 
\begin{eqnarray}\label{4.13}
\delta{J}_\gamma(d)= \left( -\int\limits_0^T q_{x}   v_{x} dt,\delta d \right)+\gamma \Big(d, \delta d\Big)_{H^1(\Omega)}+O\left(\|\delta d\|_{H^1(\Omega)}^2\right).
\end{eqnarray}
In order to bring $H^1(\Omega)$ inner product in the first term of RHS we use the equation \eqref{4.14} satisfied by $\Lambda$. Integration by parts with respect to spatial variable shows that
\begin{eqnarray} \label{4.13b}
\left( -\int\limits_0^T q_{x}   v_{x} dt,\delta d \right) = -\int\limits_{\Omega} (\Lambda^{''} - \Lambda) \delta d dx = \Big(\Lambda , \delta d\Big)_{H^1(\Omega)}.
\end{eqnarray} 
Consequently from $\eqref{4.13}$ and $\eqref{4.13b}$, $J_{\gamma}(d)$ is Frechet differentiable and the derivative is given by $\Lambda + \gamma d$.
\hfill{$\Box$}
\section{Stability Results}
\quad First, we develop the prerequisite first-order optimality condition that the optimal solution must satisfy. This condition plays a major role in proving the stability results of the unknown dissipative parameter. Next, we analyze the stability results for the inverse problem of recovering the dissipative parameter $d(x)$ of $\eqref{dir}$ from the final time measurements. The role of the optimality condition is indispensable for the proof of the stability estimate.
\begin{Lem}
Suppose $(u,v), \ (p,q)$ be the solutions to $\eqref{dir}$ and \eqref{adj} respectively and $d^* \in \mathcal{D}$ be the solution to the optimal problem $\eqref{3.1a}$. Then the following variational inequality holds :
\begin{eqnarray}\label{5.1}
-\int\limits_0^T \int\limits_{\Omega}(k(x)-d^*(x)) q_{x} v_{x} dx dt+\gamma \Big(d^*,(k-d^*)\Big)_{H^1(\Omega)} \geq 0,
\end{eqnarray}
for any  $k\in \mathcal{D}$.
\end{Lem}
\noindent {\bf Proof:} Let $k \in \mathcal{D}, \ 0 \leq \beta \leq 1$ and choose an element $d_{\beta}(x)= d^*(x) + \beta(k(x)-d^*(x)) \in  \mathcal{D} $. Let 
 $(u_\beta,v_\beta,d_\beta)$ be the solution of $\eqref{dir}$ such that objective functional $\eqref{3.1}$ becomes 
 \begin{eqnarray}\label{5.1a}
 J_{\gamma}(d_{\beta})=\!\frac{1}{2}\!\int\limits_{\Omega}\!\! \left[(u_{\beta}(x,T;d_{\beta})-m_1(x))^2 \!+\!(v_{\beta}(x,T;d_{\beta})-m_2(x))^2 \right]dx\!+\frac{\gamma}{2}\| d_{\beta}\|^2_{H^1(\Omega)}.
\end{eqnarray}
By Theorem $\eqref{T4.1}$, we knew that the functional $\eqref{5.1a}$ is Fr{$\acute{\textnormal{e}}$}chet differentiable. So, we have
 \begin{eqnarray}\label{5.2}
\left.\frac{d}{d\beta}J_{\gamma}(d_{\beta}) \right|_{\beta=0}\!\!\!&=&\!\!\!\int\limits_{\Omega}\!\! {(u_{\beta}(x,T;d_{\beta})-m_1(x))\!\left.\frac{\partial {u_{\beta}}}{\partial\beta}\right|_{\beta=0} dx}+\int\limits_{\Omega}{(v_{\beta}(x,T;d_{\beta})-m_2(x))\left.\frac{\partial {v_{\beta}}}{\partial\beta}\right|_{\beta=0} dx}\nonumber\\ 
&&+ \ {\gamma}\int\limits_{\Omega}\left[d_{\beta}(x)(k-d^*)(x)+ d'_{\beta}(x)(k-d^*)'(x)\right] dx.
 \end{eqnarray}
 Taking $\eta(x,t)=\displaystyle \frac{\partial {u_\beta}}{\partial\beta}\left. \right|_{\beta=0}$, $\theta(x,t)= \displaystyle \frac{\partial {v_\beta}}{\partial\beta}|_{\beta=0}$ and $u=u_{\beta}|_{\beta=0},v=v_{\beta}|_{\beta=0}$, then the pair $(\eta,\theta)$ satisfies the following system
\begin{eqnarray} \label{5.4}
 \hspace{-.2in}\left\{ \begin{array}{ll}
 \eta_{t}+ \eta_{xxxx} =\theta, \ (x,t) \  \in \ \Omega_T,  \\[2mm]
 \theta_{t}- (d^*(x) \theta_{x})_{x}=((k(x)-d^*(x)) v_{x})_{x}, \ (x,t) \  \in \ \Omega_T,  \\[2mm]
 \eta(x,0)=0,\  \theta(x,0)=0, \ x \in \Omega, \\[2mm]
 \eta(x,t)= \eta_{xx}(x,t)=0,\ \theta(x,t)=0, \ (x,t) \in \partial \Omega_T.
 \end{array}
 \right.
 \end{eqnarray}
 Since $d^*$ is the optimal solution, we have
\begin{eqnarray}\label{5.5}
 \left.\frac{d}{d\beta}J_{\gamma}(d^* + \beta(k-d^*))\right|_{\beta=0}\geq 0, \ \ \mbox{for any} \ k \in \mathcal{D}.
 \end{eqnarray}
 From $\eqref{5.2}$, we get
 \begin{eqnarray}\label{5.6}
 \lefteqn{\hspace{-1.0in}\int\limits_{\Omega}\!(u(x,T;d^*)-m_1(x))\eta(x,T)dx+\int\limits_{\Omega}(v(x, T;d^*)-m_2(x))\theta(x,T)dx}\nonumber\\
&& \hspace{2.1in}+ \gamma \Big(d,(k-d)\Big)_{H^1(\Omega)} \geq 0, 
 \end{eqnarray}
\mbox{for any} $\ k \in \mathcal{D}$.
 Multiplying $\eqref{5.4}_a$ and $\eqref{5.4}_b$ respectively by $p$ and $q$ the solution of the adjoint system \eqref{adj}, integrating over $\Omega_T$ and applying integration by parts several times, we get
 \begin{eqnarray}
\lefteqn {\hspace{-1in}\int\limits_{\Omega}(u(x,T;d^*)-m_1(x))\eta(x,T) dx+\int\limits_{\Omega} (v(x,T;d^*)-m_2(x))\theta(x,T)dx}\nonumber \\ \label{5.7}&& \hspace{2in}=-\int\limits_{0}^T\!\!\int\limits_{\Omega}(k-d^*)(x) v_{x} q_{x} dx dt.
 \end{eqnarray}
 Substituting $\eqref{5.7}$ in $\eqref{5.6}$, we will arrive at \eqref{5.1}.
 \hfill{$\Box$} \\
 Let $(u,v,d)$ and $(\widetilde{u},\widetilde{v},\widetilde{d})$ be the two solutions of $\eqref{dir}$ and take $U=u-\widetilde{u}, \ V=v-\widetilde{v}$ and $D=d-\widetilde{d}$ satisfying the following initial boundary value problem
\begin{eqnarray} \label{6.1}
 \hspace{-.2in}\left\{ \begin{array}{ll}
 U_{t}+ U_{xxxx}=V, \ (x,t) \  \in \ \Omega_T,  \\[2mm]
 V_{t}- (\tilde{d}(x) V_{x})_{x}=(D(x) v_{x})_{x}, \ (x,t) \  \in \ \Omega_T,  \\[2mm]
 U(x,0)=0,\  V(x,0)=0, \ x \in \Omega, \\[2mm]
 U(x,t)=0,\ U_{xx}(x,t)=0,\ V(x,t)=0, \ (x,t) \in \partial\Omega_T.
 \end{array}
 \right.\end{eqnarray}
Similarly, take $P=p-\widetilde{p} \mbox{ and } Q=q-\widetilde{q} $ where $(p,q)$ and $(\widetilde{p},\widetilde{q})$ be the two adjoint solutions of $(u,v,d)$ and $(\widetilde{u},\widetilde{v},\widetilde{d})$ respectively,  satisfying
\begin{eqnarray}\label{6.2}
 \hspace{-.2in}\left\{ \begin{array}{ll}
 -P_{t}+ P_{xxxx}=0, \ (x,t) \  \in \ \Omega_T,  \\[2mm]
 -Q_{t}-(\widetilde{d}(x) Q_{x})_{x}=P+(D(x) q_{x})_{x},\ \ \ \ (x,t) \  \in \ \Omega_T,  \\[2mm]
 P(x,T)=P_T(x),\ x \in \Omega, \\[2mm]
 Q(x,T)=Q_T(x), \ x \in \Omega, \\[2mm]
 P(x,t)=0,\ P_{xx}(x,t)=0,\ Q(x,t)=0, \ (x,t) \in \partial\Omega_T,
 \end{array}
 \right.\end{eqnarray}
 where $P_T(x)=\delta u(x,T)- (m_1(x)-\widetilde{m}_1(x))$ and $Q_T(x)=\delta v(x,T)-(m_2(x)-\widetilde{m}_2(x))$.\\
 Next, we derive some priori estimates from the equations $\eqref{6.1}, \ \eqref{6.2}$ which are essential to prove the stability result.
\begin{Lem}
Assume $d, \ \widetilde{d} \in \mathcal{D}$ and $U, \ V$ be the solution of $\eqref{6.1}$. Then the following estimate holds:
\begin{eqnarray}\label{6.3}
\lefteqn{ \hspace{-0.85in}\| U(t)\|^2+\| V (t)\|^2+2\int\limits_0^T \| U_{xx}(\tau)\|^2 d\tau +\alpha_0\int\limits_0^T \| V_{x}(\tau)\|^2 dt }\nonumber\\
&&\hspace{1.5in}\leq C_4(\alpha_0 , T)(\|u_0\|^2+\|v_0\|^2)\|D\|_{L^{\infty}(\Omega)}^2,
 \end{eqnarray}
where $C_4(\alpha_0 , T)=\displaystyle \frac{(1+Te^T)}{\alpha_0}C_1(T,\alpha_0).$ 
\end{Lem}
\noindent {\bf Proof:} Multiplying $\eqref{6.1}_a$ and $\eqref{6.1}_b$ respectively by $U$ and $V$, integrating over $\Omega$ and then applying integration by parts several times, we get
\begin{eqnarray*}
\frac{1}{2}\frac{d}{dt}\int\limits_{\Omega}(|U|^2 +|V|^2)dx+ \int\limits_{\Omega} | U_{xx}|^2 dx +\int\limits_{\Omega} \tilde{d}| V_{x}|^2 dx= \int\limits_{\Omega} UV dx-\int\limits_{\Omega} D(x) v_{x} V_{x} dx. \\[1mm]
\end{eqnarray*}
Applying Cauchy's inequality and H\"older's inequality, we get
\begin{eqnarray}\label{6.4}
\frac{d}{dt}\left(\|U(t)\|^2\!+\|V(t)\|^2 \right) +2\|  U_{xx}(t)\|^2 + \alpha_0 \| V_{x}(t)\|^2 \nonumber \\ \leq \left(\|U(t)\|^2\!+\|V(t)\|^2 \right) + \frac{1}{\alpha_0} \|D\|^2_{L^{\infty}(\Omega)} \|v_x(t)\|^2.
\end{eqnarray}
Using Gronwall's inequality and substituting \eqref{2.1} gives
\begin{eqnarray}\label{6.5}
\|U(t)\|^2 +\|V(t)\|^2 \leq \frac{e^T}{\alpha_0} C_1(T,\alpha_0) (\| u_{0}\|^2 +\|v_{0}\|^2) \|D\|^2_{L^\infty(\Omega)}.
\end{eqnarray}
Integrating \eqref{6.4} over $(0,T)$ and substituting \eqref{6.5}, we can attain \eqref{6.3}. 
\hfill{$\Box$}
\begin{Lem}
Assume $d, \ \widetilde{d} \in \mathcal{D}, \ (m_1,m_2) \in L^2(\Omega) \times L^2(\Omega), \ (\widetilde{m}_1, \widetilde{m}_2) \in L^2(\Omega) \times L^2(\Omega)$ and $P, \ Q$ be the solution of $\eqref{6.2}$. Then the following estimate holds:
\begin{eqnarray}\label{6.6}
\lefteqn{\hspace{-0.55in}\|P(t)\|^2\ +\|Q(t)\|^2+2\int\limits_0^T \| P_{xx}(\tau)\|^2 d\tau + \alpha_0 \int\limits_0^T \| Q_{x}(\tau) \|^2 d\tau }\nonumber\\ 
  &&\hspace{-0.2in} \leq (1+Te^T)(\|P_T\|^2+\|Q_T\|^2) + \frac{1}{2\alpha_0^2}(1+Te^T)^2(\|p_T\|^2+\|q_T\|^2)\|D\|^2_{L^{\infty}(\Omega)}. 
\end{eqnarray}
\end{Lem}
\noindent {\bf Proof:} Multiplying $\eqref{6.2}_a$ and $\eqref{6.2}_b$ respectively by $P$ and $Q$ and integrating over $\Omega$, we get 
\begin{eqnarray}\label{6.7}
-\frac{d}{dt}\left(\|P(t)\|^2+\|Q(t)\|^2\right)+2\| P_{xx}(t)\|^2+\alpha_0\| Q_{x}(t)\|^2 \nonumber \\ \leq \left(\|P(t)\|^2+\|Q(t)\|^2 \right) + \frac{1}{\alpha_0} \|D\|^2_{L^{\infty}(\Omega)} \| q_{x}(t)\|^2.
\end{eqnarray}
 Applying Gronwall's inequality over $(t,T)$ and using $\eqref{4.8a}$, we have
\begin{eqnarray}\label{6.8}
\lefteqn{\hspace{-3.0in}\|P(t)\|^2 +\|Q(t)\|^2 \leq  e^T \!\! \left(\|P_T\|^2+\|Q_T\|^2+\frac{1}{2\alpha_0^2}(1+Te^T)(\|p_T\|^2+\|q_T\|^2)\|D\|^2_{L^{\infty}(\Omega)}\!\! \right). }
\end{eqnarray}
Integrating \eqref{6.7} over $(t,T)$ and substituting \eqref{6.8}, we arrive at \eqref{6.6}.\hfill{$\Box$}
\begin{Rem}\label{rem1}
(see, \cite{sal}, page no.422) Since the embedding $H^1(a,b) \hookrightarrow L^{\infty}(a,b)$ is continuous, there exists a constant $C^*(a,b)$ such that 
$$\|u\|_{L^{\infty}(a,b)} \leq C^*(a,b)\|u\|_{H^1(a,b)},$$
where $C^*(a,b)=\displaystyle \sqrt{2} \max \{(b-a)^{-\frac{1}{2}},(b-a)^{\frac{1}{2}} \}$. As we have $a=0, \ b=1$ it shows that
$$\|u\|_{L^{\infty}(0,1)} \leq \sqrt{2} \|u\|_{H^1(0,1)}.$$
\end{Rem}
\begin{thm}
Let $d,\widetilde{d} \in \mathcal{D}$ be the minimizer of the functional $J_{\gamma}(\cdot)$ defined in \eqref{3.1} corresponding to the measured outputs $(m_1,m_2) \in L^2(\Omega) \times L^2(\Omega)$ and $(\widetilde{m}_1,\widetilde{m}_2) \in L^2(\Omega) \times L^2(\Omega)$. Then there exist a time $T_0$ and a constant $L(T_0)$ such that, for $T\geq T_0,$ the following stability holds :
\begin{eqnarray}\label{6.9}
\|d-\widetilde{d}\|_{H^1(\Omega)}\leq L(T_0)\left(\|m_1-\widetilde{m}_1\|+\|m_2-\widetilde{m}_2\|\right),
\end{eqnarray}
where $L(T_0)=\displaystyle \left(\frac{2(1+Te^T)}{\gamma \alpha_0 }\right)^{\frac{1}{2}}.$
\end{thm}
\noindent {\bf Proof:} Let $(u,v), \ (\widetilde{u},\widetilde{v})$ be the solutions of $\eqref{dir}$ and $(p,q), \ (\widetilde{p},\widetilde{q})$ be the adjoint solutions of the system $\eqref{adj}$ with the corresponding coefficients $d, \ \widetilde{d}$. Replace $k(x) \mbox{ by } \widetilde{d}(x)$ in \eqref{5.1} to get
\begin{eqnarray}\label{6.10}
-\int\limits_{0}^T\int\limits_{\Omega}(\widetilde{d}(x)-d(x)) v_{x} q_{x} dx dt +\gamma (d(x),(\widetilde{d}(x)-d(x)))_{H^1(\Omega)} &\geq &0.
\end{eqnarray}
\mbox{By changing $k(x) \mbox{ by } d(x)\mbox{ and } d(x) \mbox{ by } \widetilde{d}(x)$, we have}
\begin{eqnarray}\label{6.11}
-\int\limits_{0}^T \int\limits_{\Omega}(d(x)-\widetilde{d}(x)) \widetilde{v}_{x} \widetilde{q}_{x} dx dt +\gamma (\widetilde{d}(x),(d(x)-\widetilde{d}(x)))_{H^1(\Omega)} &\geq &0. 
\end{eqnarray}
Adding $\eqref{6.10} \mbox{ and } \eqref{6.11}$, we can write
\begin{eqnarray}\label{6.12}
\gamma \|d-\widetilde{d}\|^2_{H^1(\Omega)} \leq \int\limits_{0}^T\int\limits_{\Omega}D(x) V_{x} q_{x} dx dt+\int\limits_{0}^T\int\limits_{\Omega}D(x)\widetilde{v}_{x} Q_{x} dx dt:=I_3+I_4.
\end{eqnarray}
On applying Cauchy's inequality and H\"older's inequality to $I_3$, we get
\begin{eqnarray*}
I_3\leq \frac{1}{2}\|D\|^2_{L^{\infty}(\Omega)}\| q_{x}(t)\|^2_{L^2(0,T;L^2(\Omega))} + \frac{1}{2}\| V_{x}\|^2_{L^2(0,T;L^2(\Omega))}.
\end{eqnarray*}
Substituting \eqref{4.8a} and \eqref{6.3}, we obtain
\begin{eqnarray}\label{6.13}
I_3 \leq \frac{(1+Te^T)}{4\alpha_0}(\|p_T\|^2+\|q_T\|^2)\|D\|^2_{L^{\infty}(\Omega)}+\frac{C_4(T,\alpha_0)}{\alpha_0}(\|u_0\|^2+\|v_0\|^2)\|D\|^2_{L^{\infty}(\Omega)}.
\end{eqnarray}
Similarly, using \eqref{2.1} and \eqref{6.6}, we get
\begin{eqnarray}\label{6.14}
 I_4 &\leq& \frac{C_1(T,\alpha_0)}{2} (\|u_0\|^2+\|v_0\|^2) \|D\|^2_{L^{\infty}(\Omega)}+ \frac{(1+Te^T)}{2\alpha_0}(\|P_T\|^2+\|Q_T\|^2)\nonumber\\
 &&+\frac{(1+Te^T)^2}{4\alpha_0^3}(\|p_T\|^2+\|q_T\|^2)\|D\|^2_{L^{\infty}(\Omega)}.
\end{eqnarray}
Further, by applying $\eqref{4.3}$, we obtain
\begin{eqnarray*}
\|P_T\|^2+\|Q_T\|^2\!\!\! &\leq&\!\!\! 2 \left(\|\delta u(\cdot,T)\|^2+\|\delta v(\cdot,T)\|^2+\|m_1-\widetilde{m}_1\|^2+\|m_2-\widetilde{m}_2\|^2 \right) \\[2mm]
\!\!\!&\leq&\!\!\! 2\left(\frac{e^T(1+Te^T)}{2\alpha_0^2}(\|u_0\|^2+\|v_0\|^2)\|\delta d\|^2_{L^{\infty}(\Omega)}+\|m_1-\widetilde{m}_1\|^2+\|m_2-\widetilde{m}_2\|^2 \right),\\[2mm]
&&\hspace{-1.4in}\mbox{ and by \eqref{2.6}, one obtains that}\\[2mm]
&&\hspace{-1.4in}\|p_T\|^2+\|q_T\|^2 \leq 2\left(2e^T\|u_{0}\|^2+2e^T\|v_{0}\|^2+\|m_1\|^2+\|m_2\|^2\right).
\end{eqnarray*}
Substituting $\eqref{6.13} \mbox{ and } \eqref{6.14}$ in $\eqref{6.12}$ and then using Remark \ref{rem1}, we attain
\begin{eqnarray}\label{6.15}
 \|d-\widetilde{d}\|^2_{H^1(\Omega)} &\leq& \frac{\beta_1(T)}{\gamma}\|d-\widetilde{d}\|^2_{H^1(\Omega)} +\left(\frac{1+Te^T}{\gamma \alpha_0}\right)(\|m_1-\widetilde{m}_1\|^2+\|m_2-\widetilde{m}_2\|^2),
\end{eqnarray}
\begin{eqnarray}\label{6.15a}
\mbox{where }\beta_1(T)&=&\displaystyle \frac{1}{2\alpha_0^3}\left[(1+Te^T)(2\alpha_0^2+1+Te^T)(2e^T\|u_0\|^2+2e^T\|v_0\|^2+\|m_1\|^2+\|m_2\|^2) \right. \nonumber\\[2mm]  && \left. +((2\alpha_0(1+Te^T)+\alpha_0^2)(1+Te^T+4CTe^T)(1+\alpha_0)+2e^T(1+Te^T)^2) \right. \nonumber\\[2mm] &&  \left. \times (\|u_0\|^2+\|v_0\|^2) \right].
\end{eqnarray}
Choose $T_0$ such that $\displaystyle \frac{\beta_1(T_0)}{\gamma}\leq \frac{1}{2}$, we obtain 
\begin{eqnarray}\label{6.16}
\|d-\widetilde{d}\|^2_{H^1(\Omega)} \leq \frac{2(1+T_0e^{T_0})}{\gamma \alpha_0}(\|m_1-\widetilde{m}_1\|^2+\|m_2-\widetilde{m}_2\|^2),
\end{eqnarray}
which gives \eqref{6.9} and concludes the stability result.
\hfill{$\Box$}
\par In order to validate the stability result \eqref{6.9} and condition \eqref{6.15a}, we determine the upper bound for $T$ such that the inequality \eqref{6.9} holds true. From the well known inequalities for $T>0$ and $C>0,$ $1+Te^T \geq T, \ (1+Te^T)^2 \geq T, \ 2(1+Te^T)e^T \geq T, \ 2(1+Te^T)^2e^T \geq T, \ 1+Te^T+4CTe^T \geq T,$ and $(1+Te^T)(1+Te^T+4CTe^T) \geq T$, we deduce that 
\begin{eqnarray}
T \leq T^\ast
\end{eqnarray}
where 
$$T^\ast:=\frac{\gamma \alpha_0^3}{(\alpha_0^3+5\alpha_0^2+2\alpha_0+2)(\|u_0\|^2+\|v_0\|^2)+(2\alpha_0^2+1)(\|m_1\|^2+\|m_2\|^2)}.$$
By varying the values of the regularization parameter $\gamma,$ we analyse the stability constant $L$, for the following two examples.\\  
{\bf Example 1:} In the first numerical experiment, we take 
\begin{eqnarray*}
\hspace{0.4in} {d(x)=x(1-x) \mbox{ with }u(x,0)=0, \ v(x,0)=e^{-x}\sin(\pi x).}
\end{eqnarray*}
{\bf Example 2:} In the second numerical experiment, we take 
\begin{eqnarray*}
 \lefteqn{\hspace{-1.82in} d(x)=\left\{\!\!\begin{array}{c@{\ ,\quad}l}
0.5 & x\in (0.05,0.45]\\ \displaystyle{0.1+\frac{(\sin 2 \pi x)^2}{2}} & x\in[0.5,1]\\
0.1 & \mbox{otherwise}
\end{array}\right. }
\end{eqnarray*}
\hspace{1.2in} with $u(x,0)=0, v(x,0)=\sin(\pi x).$
\par The value of $T^*$ is obtained for Example 1 by choosing the values $\|u_0\|^2=0, \\ \|v_0\|^2=0.1963, \ \|m_1\|^2=7.3947 \times 10^{-4}, \mbox{ and } \|m_2\|^2=7.5540$ from the numerical experiments (which is discussed in the upcoming section). We do the analysis for the fixed value $\alpha_0=10^4$. The following table shows the change in values of the stability constant $L(T_0)$ depending on the change in values of the regularization parameter $\gamma>0$.  
\begin{table}[h]
\caption{Values of $T^*$ and its corresponding stability constant $L$ for \hspace{3in} Example 1}
\centering 
\begin{tabular}{|c|c|c|} 
\hline
$\gamma$ & $T^*$ & $L$ \\
\hline
$10^{-5}$ & $5.0528 \times 10^{-5}$ & $4.4722$ \\
$10^{-6}$ & $5.0528 \times 10^{-6}$ & $14.1422$ \\
$10^{-7}$ & $5.0528 \times 10^{-7}$ & $44.7214$ \\
$10^{-8}$ & $5.0528 \times 10^{-8}$ & $141.4214$ \\
$10^{-9}$ & $5.0528 \times 10^{-9}$ & $447.2136$ \\
$10^{-10}$ & $5.0528 \times 10^{-10}$ & $1414.2$ \\
\hline
\end{tabular}
\label{table2}
\end{table} 
\par The value of $T^*$ is obtained for Example 2 by choosing the values $\|u_0\|^2=0, \\ \|v_0\|^2=0.5, \ \|m_1\|^2=4.0266 \times 10^{-5}, \mbox{ and } \|m_2\|^2=0.3576$ from the numerical experiments (which is discussed in the upcoming section). We do the analysis for the fixed value $\alpha_0=10^4$. The following table shows the change in values of the stability constant $L(T_0)$ depending on the change in values of the regularization parameter $\gamma>0$. 
\begin{table}[h]
\caption{Values of $T^*$ and its corresponding stability constant $L$ for \hspace{3in} Example 2}
\centering 
\begin{tabular}{|c|c|c|} 
\hline
$\gamma$ & $T^*$ & $L$ \\
\hline
$10^{-5}$ & $1.9987 \times 10^{-5}$ & $4.4722$ \\
$10^{-6}$ & $1.9987 \times 10^{-6}$ & $14.1421$ \\
$10^{-7}$ & $1.9987 \times 10^{-7}$ & $44.7214$ \\
$10^{-8}$ & $1.9987 \times 10^{-8}$ & $141.4214$ \\
$10^{-9}$ & $1.9987 \times 10^{-9}$ & $447.2136$ \\
$10^{-10}$ & $1.9987 \times 10^{-10}$ & $1414.2$ \\
\hline
\end{tabular}
\label{table3}
\end{table}
\par From Tables \ref{table2} and \ref{table3}, it is clear that when the value of the regularization parameter $\gamma$ decreases below $10^{-7}$, the time $T^*$ decreases, and the stability constant value $L$ is increasing rapidly. This indicates that the stability constant $L$ holds for some specified small interval of time $T$ and also it is of the order $\gamma^{-\frac{1}{2}}$. Ultimately, to get a consistent stability constant the regularization parameter should be greater than $10^{-7}$ which is reasonable in terms of real world applications because the value of regularization parameter actually ranges between $10^{-2}$ to $10^{-4}$.
\section{Numerical Reconstruction of the dissipative parameter}
	This section provides the numerical simulation for the unknown dissipative parameter $d(x)$. Here, we have given the graphical illustration for the convergence for the dissipative parameter $d(x)$. An iterative scheme based on Conjugate Gradient Method(CGM) is used for computing $d(x)$ by minimizing the functional $\eqref{3.1}$. The procedure for the iterative process takes the following form
\begin{eqnarray}\label{7.1}
d^{n+1}(x)=d^n(x)-\beta^n q^{n}(x),
\end{eqnarray}
where $n$ is the number of iterations, $\beta^n$ is the step length/step search size, and $q^n$ is the descent direction defined by 
\begin{eqnarray}\label{7.2}
q^n=\left\{ \begin{array}{ll}
J'(d^0),& n=0\\
J'(d^n)+\mu^n q^{n-1}, & n=1,2,3,...
\end{array}\right.
\end{eqnarray}
In \eqref{7.2}, $J'(d^n)$ denotes the gradient of the functional derived in \eqref{4.15}, $\mu^n$ is the conjugation coefficient that has various expressions such as Polak-Ribiere \cite{pow}, Fletcher-Reeves \cite{FR}, Powell-Beale method \cite{pow}, etc. In this article, we have considered the Polak-Ribiere expression \cite{pow} defined as 
\begin{eqnarray}\label{7.3}
\mu^0_{PR}=0,\ \mu^n_{PR}=\frac{(J'(d^n), J'(d^n)- J'(d^{n-1}))}{\|J'(d^{n-1})\|^2},
\end{eqnarray}
and the Fletcher-Reeves expression \cite{FR} defined as
\begin{eqnarray}\label{7.4}
\mu^0_{FR}=0,\ \mu^n_{FR}=\frac{\| J'(d^{n})\|^2 }{\| J'(d^{n-1})\|^2}.
\end{eqnarray}
In both the expressions \eqref{7.3} and \eqref{7.4}, we calculate $J'(d^{n-1})$ using \eqref{4.15}. The step search size $\beta^n$ is computed using exact line search method is defined by
\begin{eqnarray*}
\beta^n =\arg \min\limits_{\beta\geq 0} J_{\gamma}(d^{n}(x)-\beta q^n(x)),
\end{eqnarray*}
which is equivalent to 
\begin{eqnarray}\label{7.5}
\left.\frac{\partial}{\partial \beta}J_{\gamma}(d^n(x)-\beta q^n(x))\right|_{\beta=\beta^n}\!\!\!\!\!=\!0.\ 
\end{eqnarray}
By $\eqref{3.1},$ we have
\begin{eqnarray}\label{7.6}
{J_{\gamma}(d^{n}-\beta q^n)}\!\!&=&\!\!\frac{1}{2}\!\int\limits_{\Omega}\! {\left(u(x,T;d^{n}\!-\beta q^n)\!-m_1(x)\right)^{2} \!dx}\!+\frac{1}{2}\!\int\limits_{\Omega}\!{\left(v(x,T;d^{n}\!-\beta q^n)\!-m_2(x)\right)^{2}\! dx}\nonumber\\
\lefteqn{\hspace{0.3in}+\frac{\gamma}{2} \|d^{n}-\beta q^n\|^{2}_{H^1(\Omega)}.}
\end{eqnarray}
Set $q^n=\delta d^n$ and by Taylor's expression, we obtain
\begin{eqnarray}\label{7.7}
\left\{\begin{array}{c}
u(x,T;d^{n}-\beta q^n)\approx u(x,T;d^{n})-\displaystyle{\frac{\partial u(x,T;d^{n})}{\partial d^n}}\beta q^n\approx u(x,T;d^{n})-\beta \delta u(x,T;d^{n}),\\[3mm]
v(x,T;d^{n}-\beta q^n)\approx v(x,T;d^{n})-\displaystyle{\frac{\partial v(x,T;d^{n})}{\partial d^n}}\beta q^n\approx v(x,T;d^{n})-\beta \delta v(x,T;d^{n}),
\end{array}
\right.
\end{eqnarray}
where $\delta u(x,T;d^{n}) \mbox{ and } \delta v(x,T;d^{n})$ represent the solution of the sensitivity problem $\eqref{4.1}$ at $t=T$. Substitute \eqref{7.7} in \eqref{7.6} and then by \eqref{7.5}, we get
\small{\begin{eqnarray}\label{7.8}
\beta^n\!=\!\frac{\int\limits_{\Omega}\![\left(u(x,T;d^{n})\!-\!m_1(x)\right)\delta u(x,T;d^{n}) \!+\!\left(v(x,T;d^{n})\!-\!m_2(x)\right)\delta v(x,T;d^{n}) ]dx\!+\!\gamma (d^n,q^n)_{H^1(\Omega)}}{\|\delta u(x,T;d^{n})\|^2+\|\delta v(x,T;d^{n})\|^2\!+\gamma\| q^n\|^2_{H^1(\Omega)}}.
\end{eqnarray}}
\subsection{Algorithm}
\begin{description}
\item[{\it Step $1$}]Choose an initial guess $d^0$ for the unknown dissipative parameter $d(x)$ and set $n=0.$
\item[{\it Step $2$}]Solve the direct problem $\eqref{dir}$ using the finite difference scheme to calculate \\$(u(x,t;d^n),v(x,t;d^n)$, and then determine the functional $J_\gamma(d^n)$ given in \eqref{3.1}.
\item[{\it Step $3$}]Solve the adjoint problem $\eqref{adj}$ to compute $(p(x,t;d^n),q(x,t;d^n)$, and the Fr\'echet gradient of the functional $J^{\prime}(d^n)$ derived in $\eqref{4.15}$. Further, we find the conjugate coefficient $\mu^n$ from \eqref{7.3} and \eqref{7.4}, and the descent direction $q^n$ in $\eqref{7.2}$.
\item[{\it Step $4$}] Solve the sensitivity problem \eqref{4.1} to determine $(\delta u(x,t;d^n),\delta v(x,t;d^n)$ by choosing $\delta d^n=q^n$, and calculate the step search size $\beta^n$ from $\eqref{7.8}$.
\item[{\it Step $5$}]Update the coefficient $d^{n+1}$ in \eqref{7.1} for $n=1,2,3,...$ Revert to {\it Step $2$} and repeat the procedure untill the following stopping criterion for the iterative procedure is satisfied.
\begin{eqnarray}\label{7.9}
J_\gamma(d^n)\leq \varrho \epsilon,
\end{eqnarray}
\end{description}
where $\displaystyle \varrho>1$ (for example 1.01) and $\epsilon$ denotes the noise added to the measured output data which is given by $$\epsilon=\frac{1}{2}(\|u_h-u_e\|^2+\|v_h-v_e\|^2).$$
All the integrals connected with the algorithm are numerically computed using Simpson's rule. The direct and adjoint problem are solved using Finite Difference Method and the value for $\epsilon$ mentioned in the discrepancy principle is provided in the following section.
\subsection{Numerical Results and Discussions} 
 \par In this section, the numerical experiments were performed for the iterative CGM algorithm which was mentioned in the previous section for the inverse problem $\eqref{dir}$-$\eqref{1.2}.$ Here, we used Finite difference scheme (\cite{ozi},\cite{sam}) to solve the direct and adjoint problems. To proceed with this scheme, we initially need to discretize the space and time interval into uniform grids that is, $\Delta x=1/N, \Delta t=T/N$. We apply forward difference scheme for the time derivative and a control volume approach for the discretization of spatial derivative of $\eqref{dir}$ which leads to the following discrete form
\normalsize{\begin{eqnarray}\label{7.10}
\displaystyle \left\{\begin{array}{ll}
\frac{u_i^{j+1}-u_i^{j}}{\Delta t}+\frac{1}{(\Delta x)^4}(u_{i+2}^{j+1}-4u_{i+1}^{j+1}+6u_{i}^{j+1}-4u_{i-1}^{j+1}+u_{i-2}^{j+1})=v_i^j,\\[2mm]
 \frac{v_i^{j+1}-v_i^{j}}{\Delta t}-\frac{1}{(\Delta x)^2}\left(\left(\frac{d(i+1)+d(i)}{2}\right)(v_{i+1}^{j+1}-v_{i}^{j+1})- \left(\frac{d(i)+d(i-1)}{2}\right)(v_{i}^{j+1}-v_{i-1}^{j+1}) \right)=0,\\[2mm]
 \end{array}\right.
 \end{eqnarray}}
 where $u_i^j:=u(x_i,t_j),\mbox{ in which } x_i=i \Delta x \mbox{ for }i=0,1,2,\cdots, N,\mbox{ and } t_j=j \Delta t \mbox{ for }j=0,1,2,\cdots, N,$ and $d_i:=d(x_i).$ 
In a similar manner, we can apply this scheme for adjoint problem $\eqref{adj}.$ 
Let $u_h, v_h$ be the noisy measured data and $u_e,  v_e $ be the exact solution of the direct problem. The noisy measured data $u_h, \ v_h$ are defined as follows
$$\displaystyle  \ u_h:=u_e+\omega p \max\limits_{[0,1]}abs(u_e),\mbox{ and } v_h:=v_e+\omega p \max\limits_{[0,1]}abs(v_e)$$
where $p$ represents the percentage of noise and $\omega$ is generated by MATLAB function\\ $randn(1,N)$ with mean $0$ and variance $1.$ In all the numerical examples, we take $T=1$, $N=100,\mbox{ the regularization parameter } \gamma=1\times 10^{-6}$.

Finally, the error at each iteration for the dissipative parameter $d(x)$ is defined as
$$Error:=\|d-d^n\|.$$
\par The numerical results for the dissipative parameter $d(x)$ given in Example 1 and 2 for percentages of noise $p=1\%,3\%,5\%$ are obtained using the CGM and are stopped in respect to the iteration numbers provided in Table \ref{tab:1} $\&$ \ref{tab:4}. We used both the Polak-Ribiere(PR) and Fletcher-Reeves(FR) conjugate coefficient formulas mentioned in $\eqref{7.3}$ and $\eqref{7.4}$ respectively to compare the results and the value of $J$ for the noisy data.
{\small{\begin{table}[h]
\caption{Performance after iterations for Example 1}
\centering 
\begin{tabular}{|l c | c c c | c c c |} 
\hline\hline
$J({\bf\Downarrow})$ & & \multicolumn{3}{c}{$J(d)$ for FR}\vline & \multicolumn{3}{c}{$J(d)$ for PR}\vline 
\\
\cline{3-8}
& $p({\bf \Rightarrow})$ & $1\%$& $3\%$ & $5\%$ & $1\%$ & $3\%$ & $5\%$ \\
\hline
  &\raisebox{-1.5pt}{J}& \raisebox{-1.5pt}{5.24 $\times 10^{-4}$} & \raisebox{-1.5pt}{4.02 $\times 10^{-3}$} & $\raisebox{-1.5pt}{0.0104}$ &\raisebox{-1.5pt}{5.07$\times 10^{-4}$} &\raisebox{-1.5pt}{4.27 $\times 10^{-3}$} & $\raisebox{-1.5pt}{0.0104}$ \\
\raisebox{8pt}{Example 1} & iterations &$20$ &$9$ & $7$& $70$ & $28$ & $20$ \\
 \hline
\end{tabular}
\label{tab:1}
\end{table}}}
{\small{\begin{table}[h]
\caption{Performance after iterations for Example 2}
\centering 
\begin{tabular}{|l c | c c c | c c c |} 
\hline\hline
$J({\bf\Downarrow})$ & & \multicolumn{3}{c}{$J(d)$ for FR}\vline & \multicolumn{3}{c}{$J(d)$ for PR}\vline 
\\
\cline{3-8}
& $p({\bf \Rightarrow})$ & $1\%$& $3\%$ & $5\%$ & $1\%$ & $3\%$ & $5\%$ \\
\hline
  &\raisebox{-1.5pt}{J}& \raisebox{-1.5pt}{1.08 $\times 10^{-3}$} & \raisebox{-1.5pt}{1.08$\times 10^{-3}$} & \raisebox{-1.5pt}{1.14$\times 10^{-3}$} &\raisebox{-1.5pt}{3.44$\times 10^{-4}$} &\raisebox{-1.5pt}{4.28$\times 10^{-4}$} & \raisebox{-1.5pt}{8.74$\times 10^{-4}$} \\
\raisebox{8pt}{Example 2} & iterations &$16380$ &$16380$ & $14788$& $16380$ & $6729$ & $297$ \\
 \hline 
\end{tabular}
\label{tab:4}
\end{table}}}
\begin{figure}[!htb]
\centering
\begin{minipage}[b]{\linewidth}
\includegraphics[height=0.3\textwidth, width=0.5\textwidth]{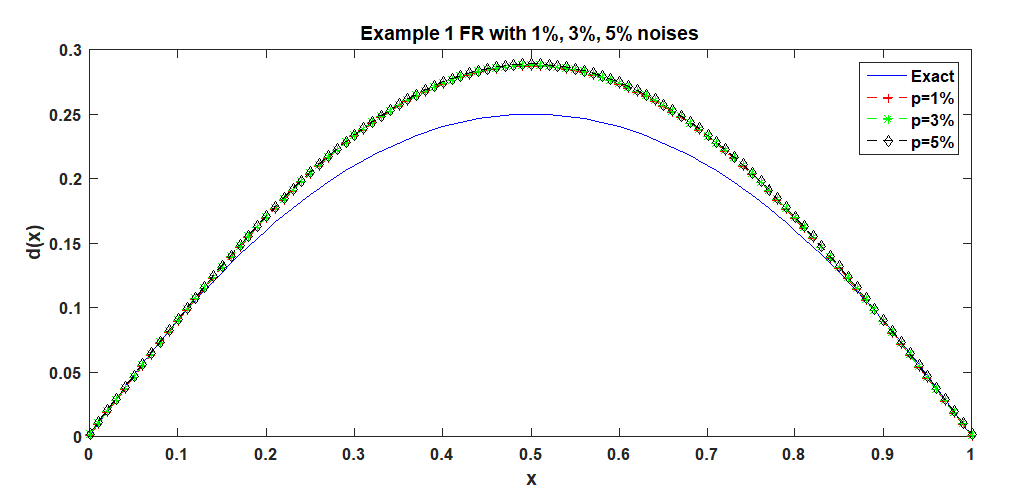}
\includegraphics[height=0.3\textwidth, width=0.5\textwidth]{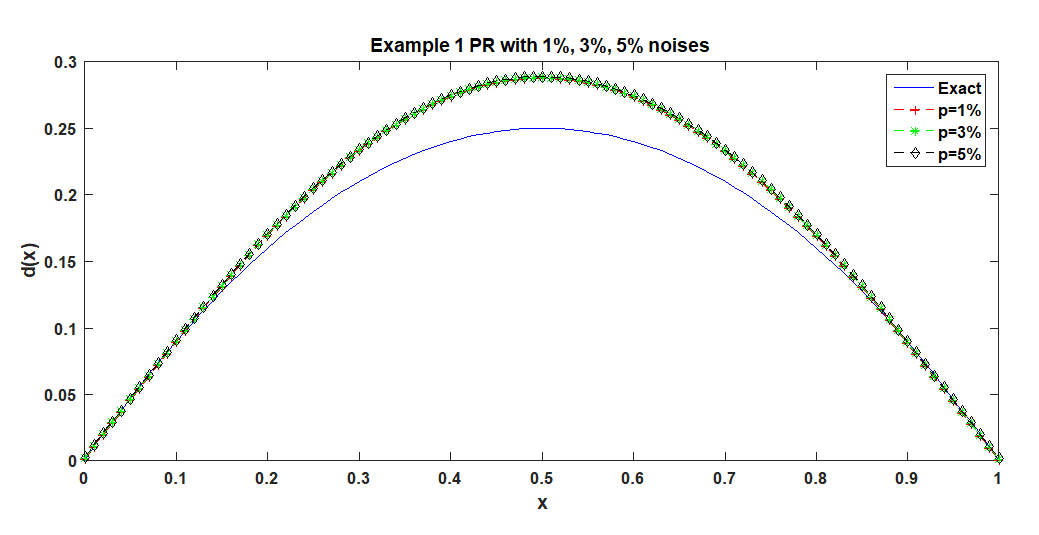}
\vspace{-0.3in}
\caption{Exact and numerical results for the reconstruction of the dissipative parameter $d(x)$ with $1\%, 3\% \ \& \ 5\%$ noisy data for Example 1}\label{figure:1}
\end{minipage}
\end{figure}
\begin{figure}[htb!]
\begin{minipage}[htb]{\textwidth}
\includegraphics[height=0.3\textwidth, width=0.5\textwidth]{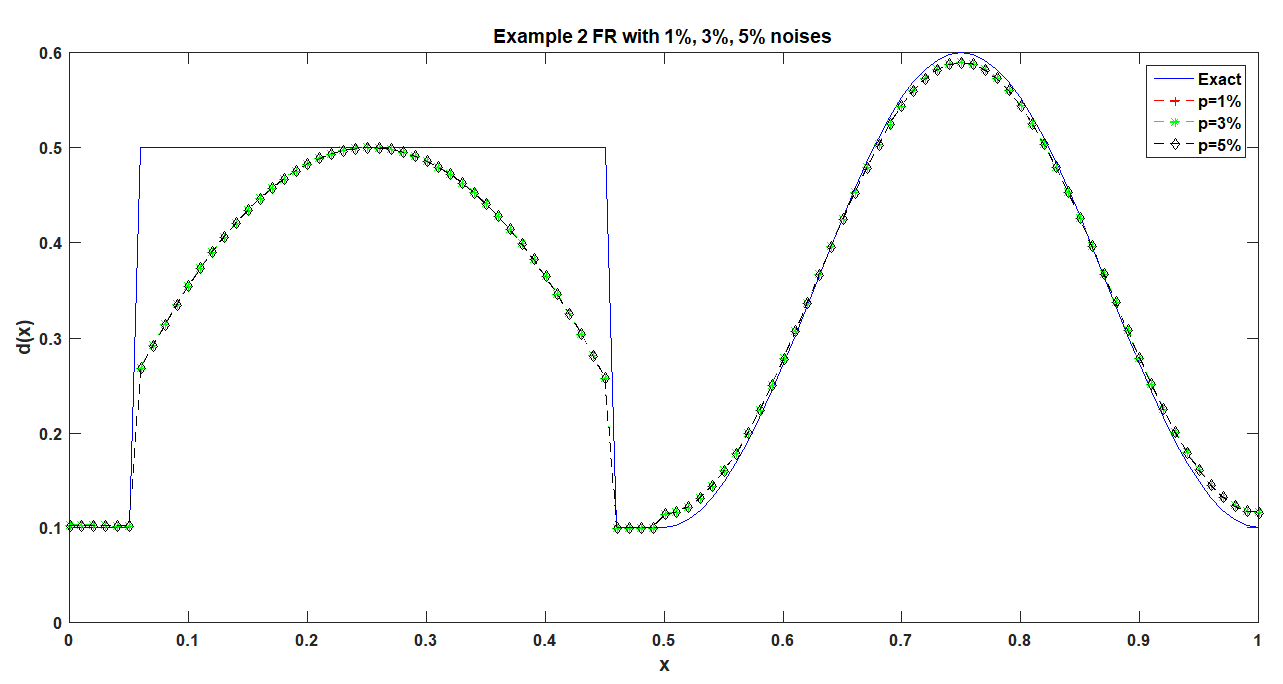} 
\includegraphics[height=0.3\textwidth, width=0.5\textwidth]{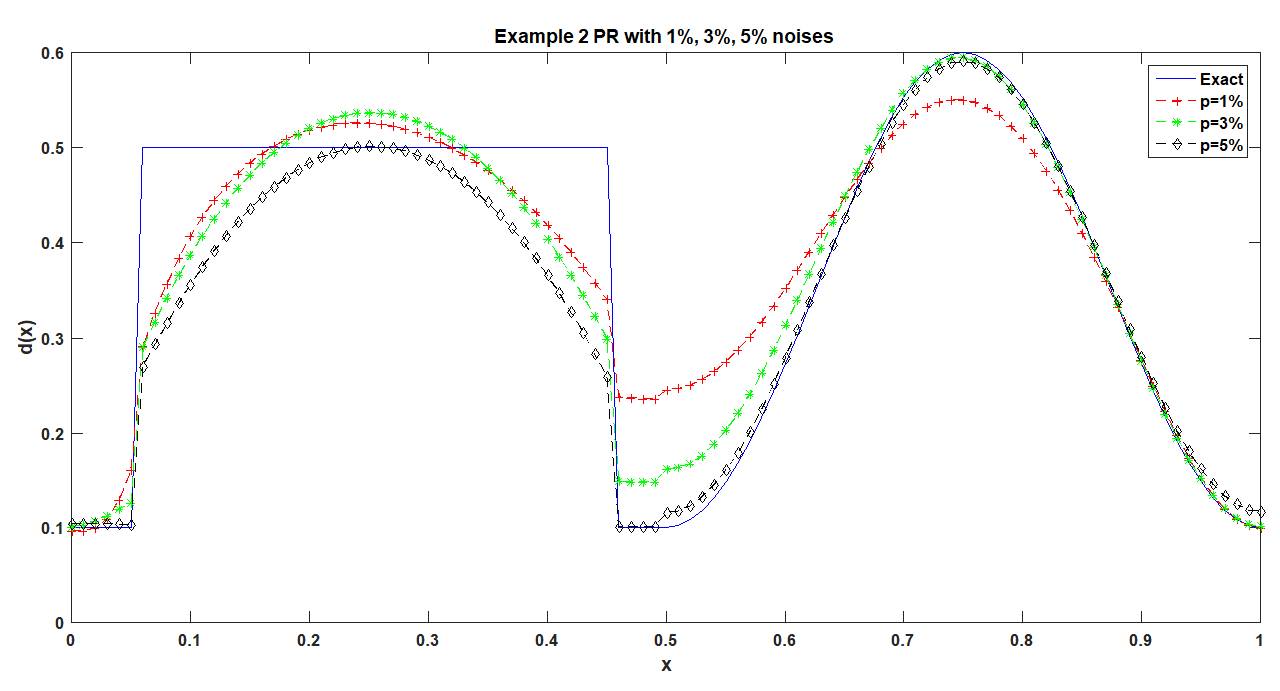}
\vspace{-0.3in}
\caption{Exact and numerical results for the reconstruction of the dissipative parameter $d(x)$ with $1\%, 3\% \ \& \ 5\%$ noisy data for Example 2}\label{figure:2}
\end{minipage}
\end{figure}
\par The numerical solutions of the unknown dissipative parameter $d(x)$ for various percentages of noise are represented in Figure-\ref{figure:1} and \ref{figure:2}. The LHS and RHS of Figure-\ref{figure:1} and \ref{figure:2} depict the comparison of unknown dissipative parameter $d(x)$ for distinct noisy data with conjugate coefficients mentioned in \eqref{7.3} and \eqref{7.4} respectively. 
\begin{figure}[htb!]
\begin{minipage}[htb]{\textwidth}
\includegraphics[height=0.3\textwidth, width=0.5\textwidth]{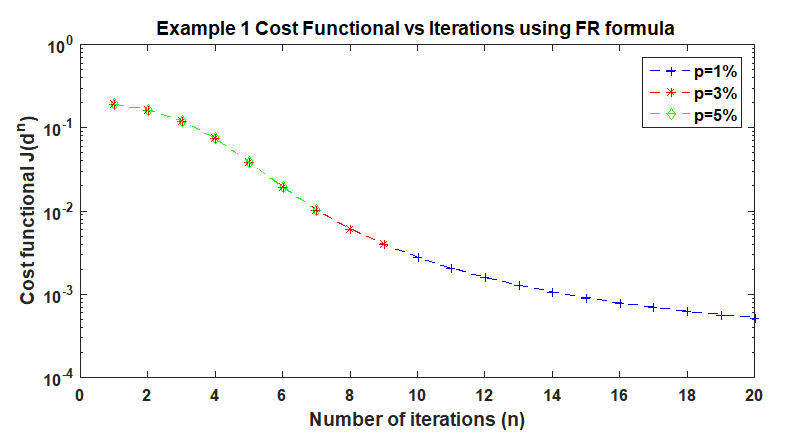} 
\includegraphics[height=0.3\textwidth, width=0.5\textwidth]{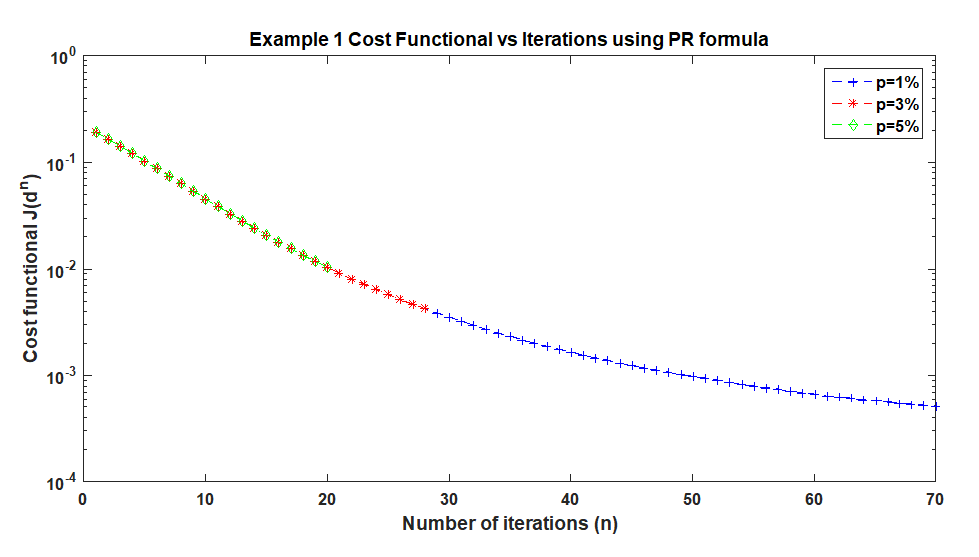}
\vspace{-0.3in}
\caption{Cost functional vs Number of iterations with $1\%, 3\% \ \& \ 5\%$ noisy data for Example 1}\label{figure:3}
\end{minipage}
\end{figure}
\begin{figure}[htb!]
\begin{minipage}[htb]{\textwidth}
\includegraphics[height=0.3\textwidth, width=0.5\textwidth]{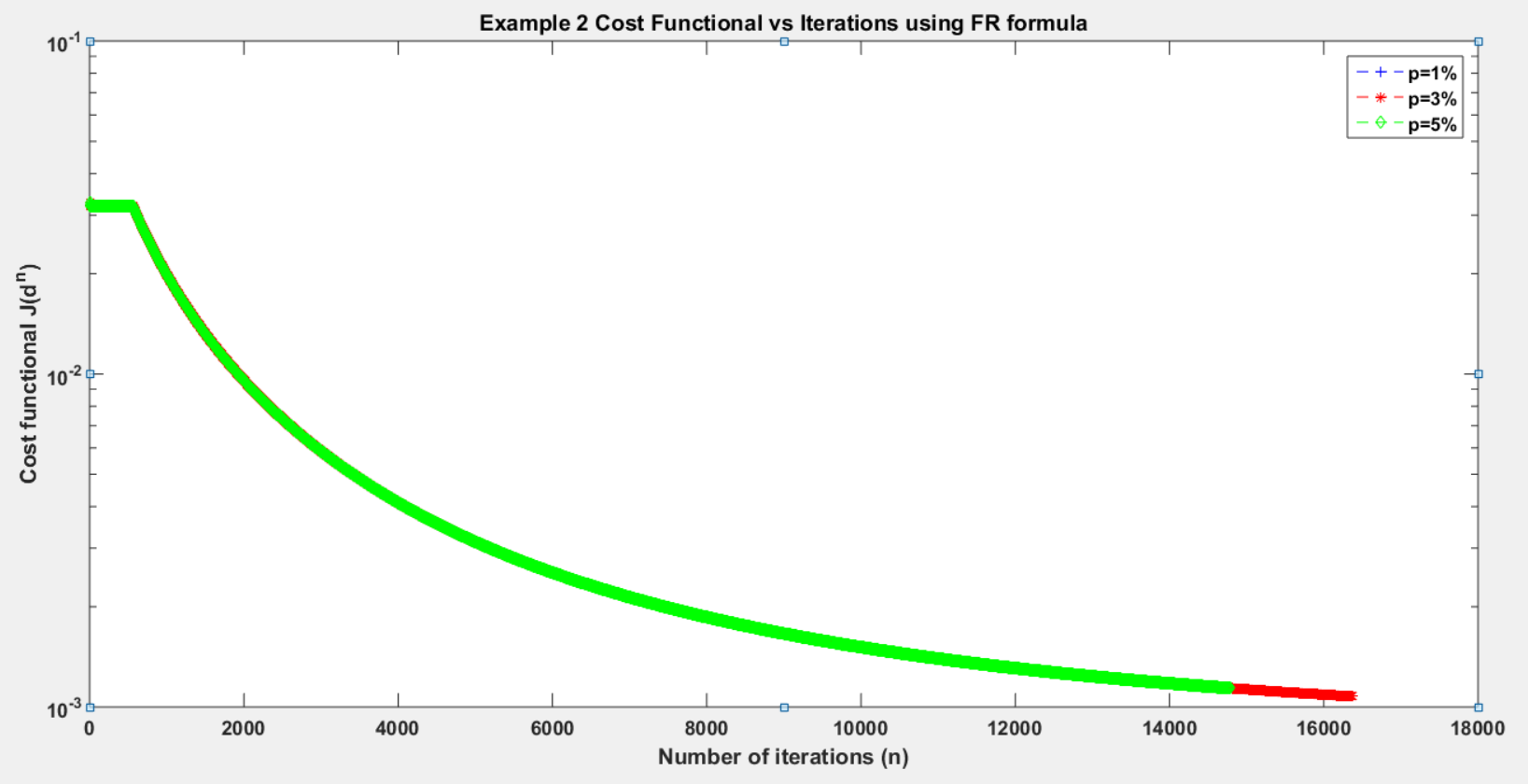} 
\includegraphics[height=0.3\textwidth, width=0.5\textwidth]{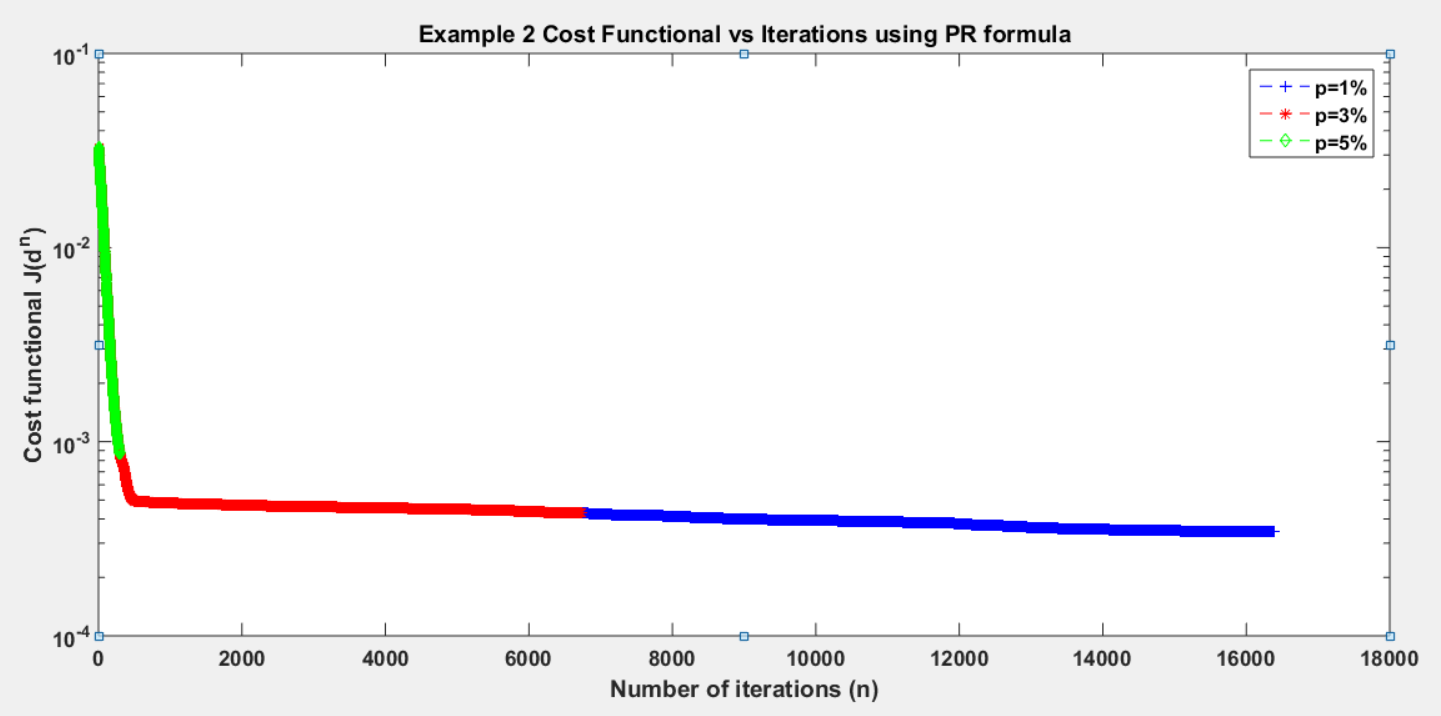}
\vspace{-0.3in}
\caption{Cost functional vs Number of iterations with $1\%, 3\% \ \& \ 5\%$ noisy data for Example 2}\label{figure:4}
\end{minipage}
\end{figure}
\begin{figure}[!htb]
\centering
\begin{minipage}[b]{\linewidth}
\includegraphics[height=0.3\textwidth, width=0.5\textwidth]{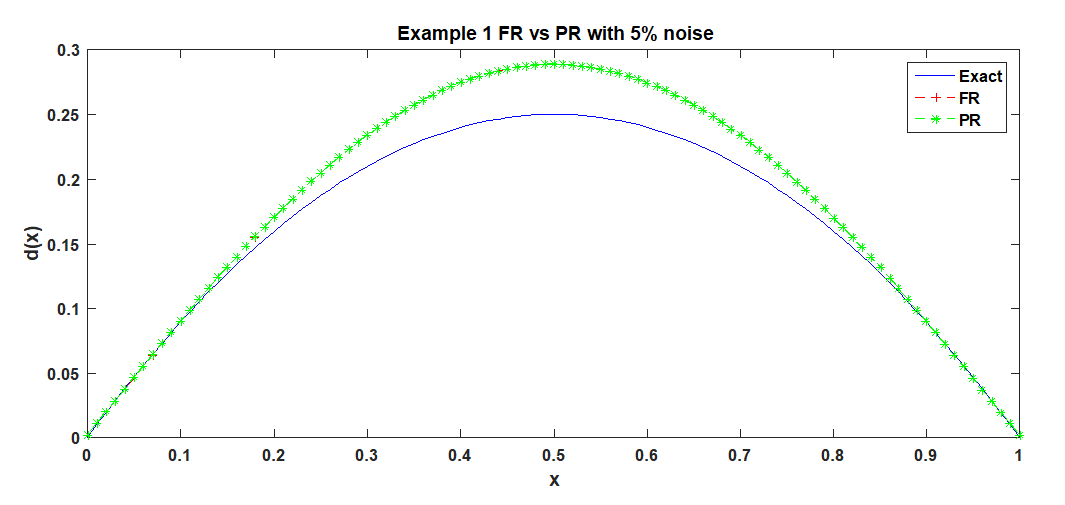}
\includegraphics[height=0.3\textwidth, width=0.5\textwidth]{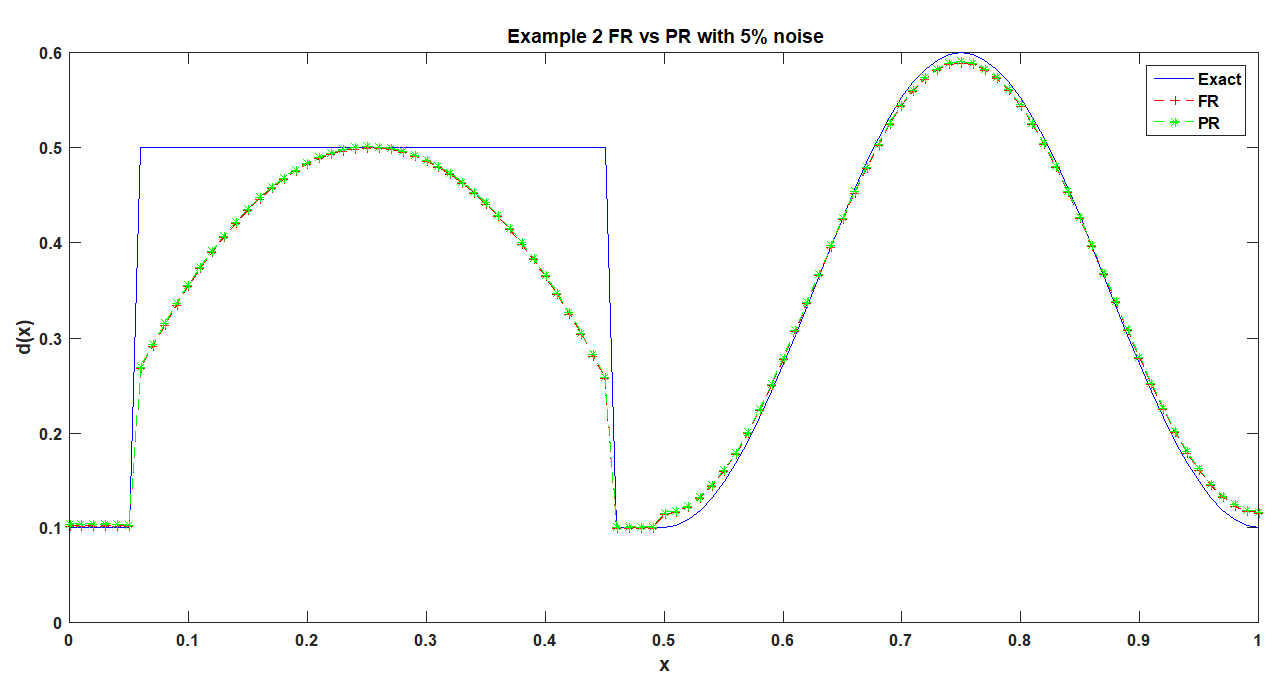}
\vspace{-0.3in}
\caption{FR vs PR for $5\%$ noisy data for the reconstruction of the dissipative parameter $d(x)$}\label{figure:5}
\end{minipage}
\end{figure}
\par Figure-\ref{figure:3} and \ref{figure:4} represent the value of the cost functional in the successive iterations for various percentages of noise using both FR formula and PR formula. As it is seen from the figure that the cost functional value gradually decreases as the number of iterations increases and it decreases further when the percentage of noise decreases. So, it is clear that the accuracy increases with decrease in noise percentage. 
\begin{figure}[!htb]
\centering
\begin{minipage}[b]{\linewidth}
\vspace{-0.5in}
\includegraphics[height=0.3\textwidth, width=0.5\textwidth]{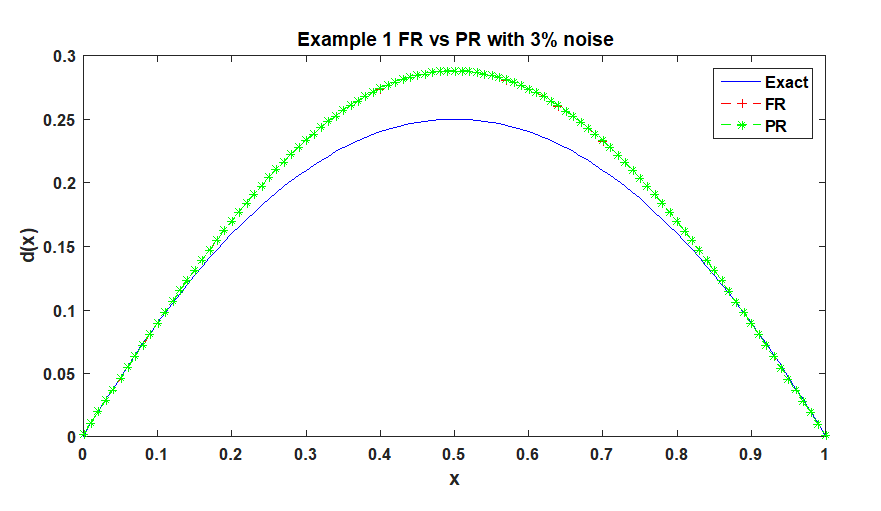}
\includegraphics[height=0.3\textwidth, width=0.5\textwidth]{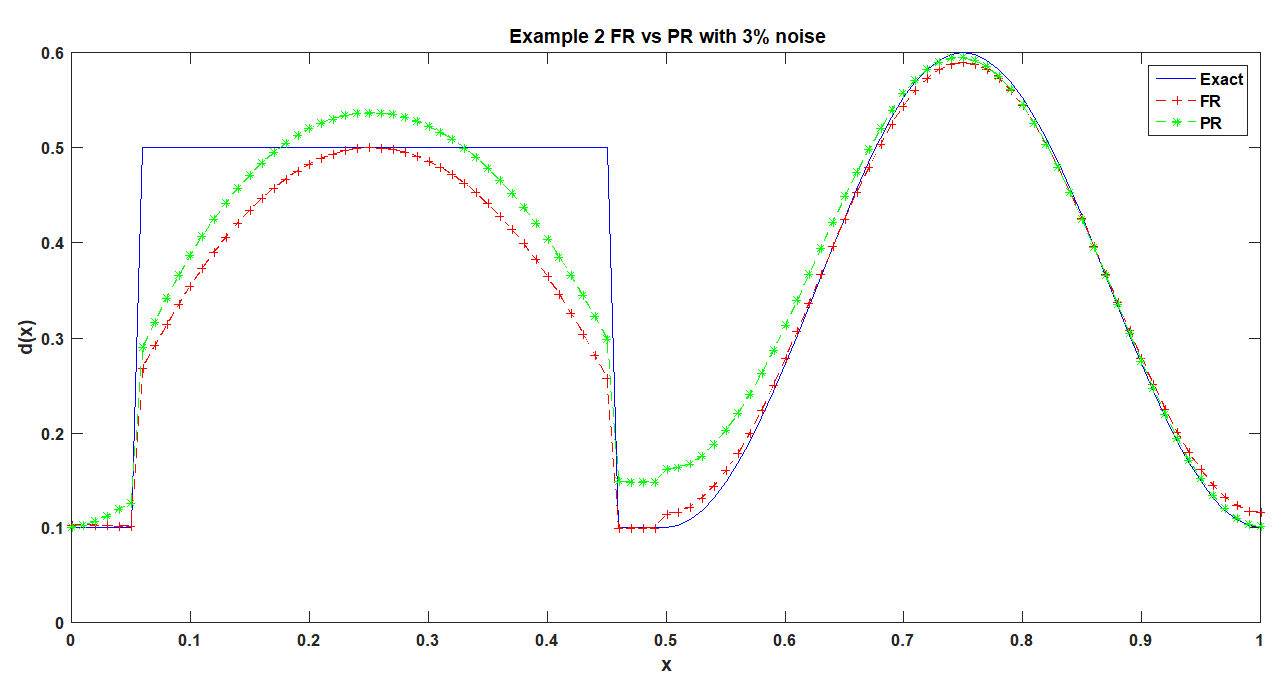}
\vspace{-0.3in}
\caption{FR vs PR for $3\%$ noisy data for the reconstruction of the dissipative parameter $d(x)$}\label{figure:6}
\end{minipage}
\end{figure}
\begin{figure}[!htb]
\centering
\begin{minipage}[b]{\linewidth}
\includegraphics[height=0.3\textwidth, width=0.5\textwidth]{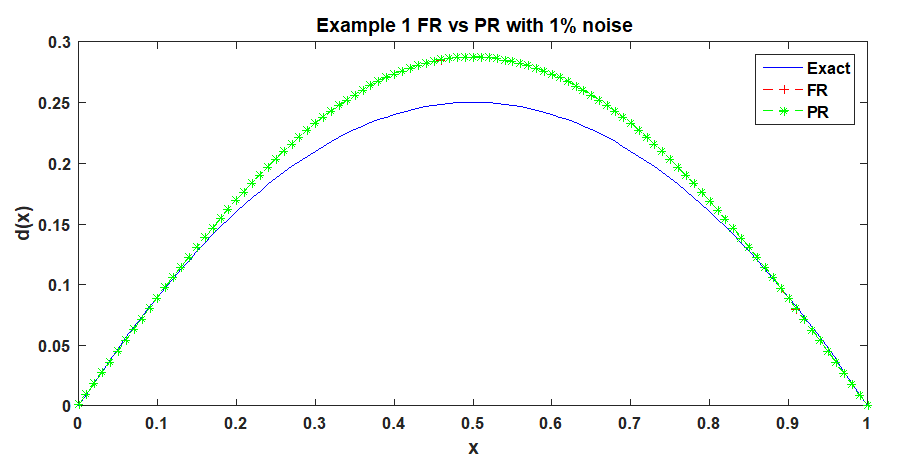}
\includegraphics[height=0.3\textwidth, width=0.5\textwidth]{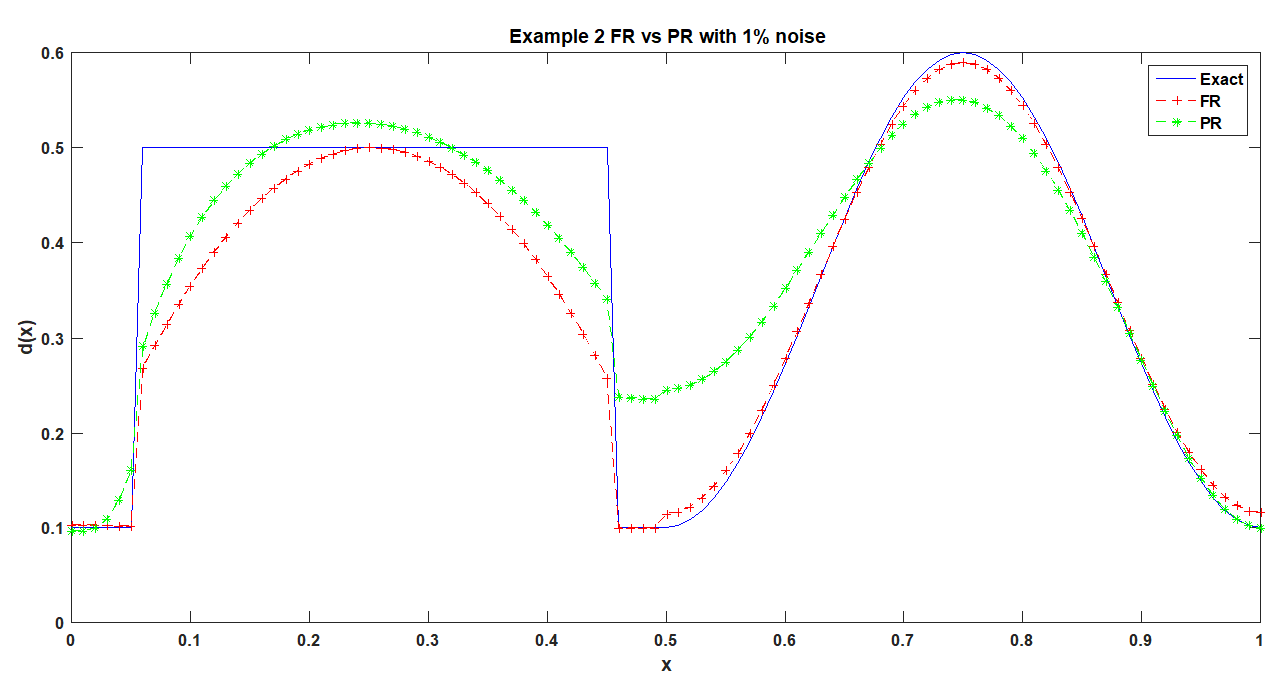}
\vspace{-0.3in}
\caption{FR vs PR for $1\%$ noisy data for the reconstruction of the dissipative parameter $d(x)$}\label{figure:7}
\end{minipage}
\end{figure}
\par Figure-\ref{figure:5}, \ref{figure:6}, \ref{figure:7} show the comparison of the Fletcher-Reeves and Polak-Ribiere formulas for different percentages of noisy data respectively for both the examples. From these figures we can see that when the noise percentages are high there is no much difference between FR formula and PR formula but when the noise is low we can clearly see that FR formula has a better accuracy than the PR formula.
\begin{figure}[!htb]
\centering
\begin{minipage}[b]{\linewidth}
\includegraphics[height=0.3\textwidth, width=0.5\textwidth]{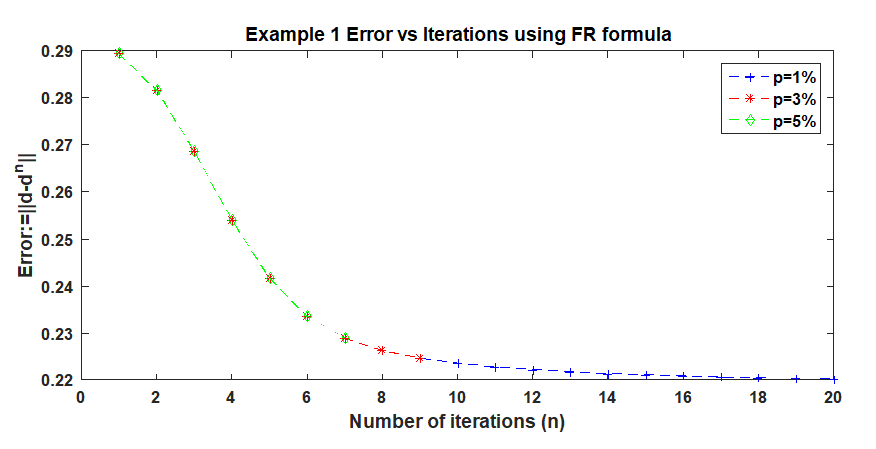}
\includegraphics[height=0.3\textwidth, width=0.5\textwidth]{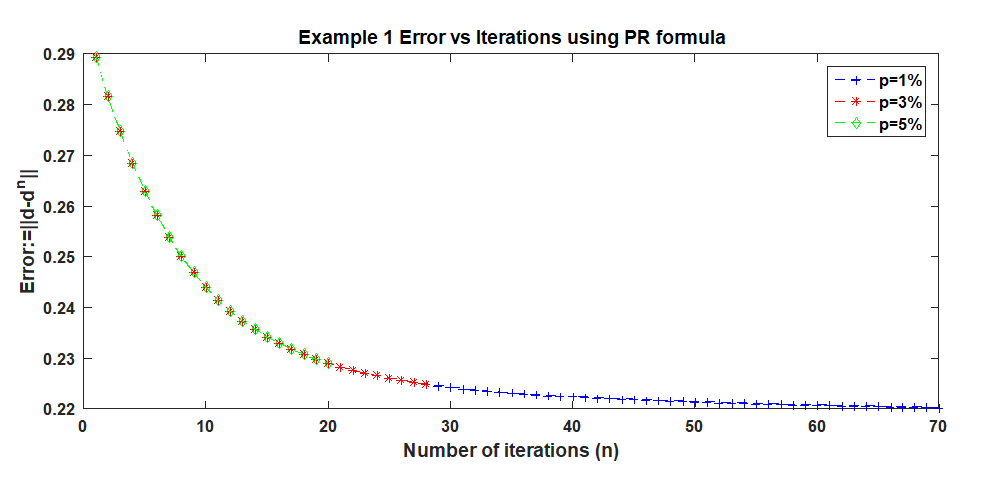}
\vspace{-0.3in}
\caption{Error vs Number of iterations for the reconstruction of the dissipative parameter $d(x)$ with $1\%, 3\% \ \& \ 5\%$ noisy data for Example 1}\label{figure:8}
\end{minipage}
\end{figure}
\begin{figure}[!htb]
\centering
\begin{minipage}[b]{\linewidth}
\includegraphics[height=0.3\textwidth, width=0.5\textwidth]{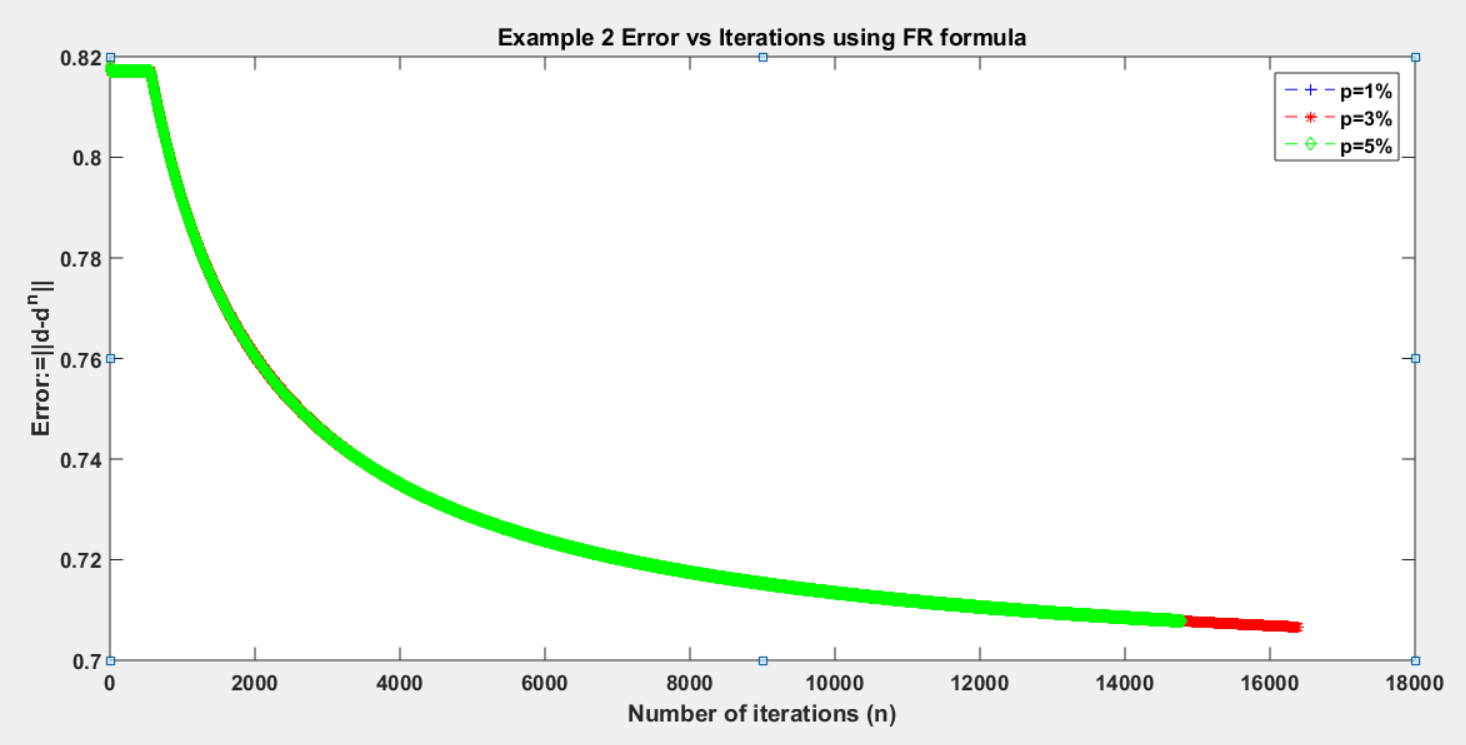}
\includegraphics[height=0.3\textwidth, width=0.5\textwidth]{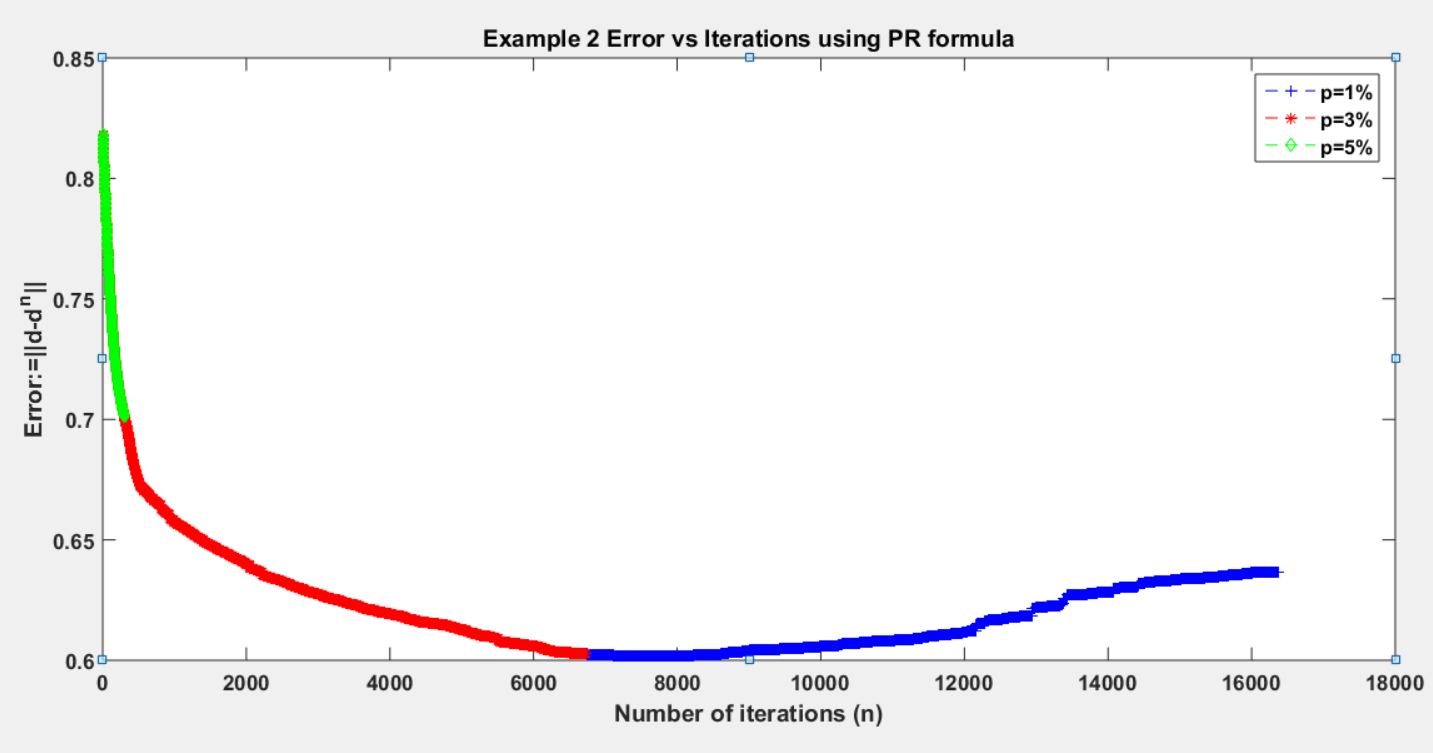}
\vspace{-0.3in}
\caption{Error vs Number of iterations for the reconstruction of the dissipative parameter $d(x)$ with $1\%, 3\% \ \& \ 5\%$ noisy data for Example 2}\label{figure:9}
\end{minipage}
\end{figure}
\par Figure-\ref{figure:8} and \ref{figure:9} give the error at each iteration for different percentages of noises using both FR formula and PR formula. From which we can notice that the error decreases as the number of iterations increases and it decreases further when the percentage of noise decreases. Except for the use of PR formula in Example 2 where for $1\%$ noisy data the error decreases and then increases and finally it will become stable. This is all because of the discontinuous nature of Example 2.
\par Finally, from the figures and Tables \ref{tab:1} $\&$\ref{tab:4} we observe that for continuous case given as Example 1 FR formula converges faster when compared with PR formula and for the discontinuous case given as Example 2 this is vice-versa. 
\section{Conclusion}
\quad This article investigates an inverse problem on recovering the dissipative parameter of a cascade system of fourth and second order PDE from final time measurement. The well-posedness of the direct and adjoint problems has been proved with the aid of Faedo-Galerkin method. The compactness of the input-output operator proves that the corresponding inverse problem is illposed. The existence of the minimizer is also shown. The explicit form of the Fr\'echet derivative of the functional is derived. The numerical treatment was considered using CGM for the Fletcher-Reeves and Polak-Ribiere conjugation coefficient formulas and the results were compared and provided in the Figures \ref{figure:5}-\ref{figure:7} and Tables \ref{tab:1}$-$\ref{tab:4}.
\section*{Acknowledgements}
The authors thank Navaneetha Krishnan Karuppusamy, Central University of Tamil Nadu, India who initiated this work and for the help made in numerical section. The first author thank the University Grants Commission, India for the financial support via CSIR- UGC Junior Research Fellowship (UGC Ref. No.: \textbf{1198/(CSIR-UGC NET DEC. 2018)}). The National Board for Higher Mathematics has funded the second author's work with research grant number\\ \textbf{02011/13/2022/R$\&$D-II/10206}. The last author was supported by SERB through the research grant No.: {\bf SR/FTP/MS-048/2011 dt. 23.06.2014.}

\end{document}